\documentclass{amsart}
\usepackage{mathtools,amssymb,graphicx,booktabs,multirow,subcaption,tikz,array}
\usetikzlibrary{arrows.meta}
\usepackage[numbers,sort&compress]{natbib}
\usepackage[hidelinks]{hyperref}
\numberwithin{equation}{section}
\numberwithin{figure}{section}
\newtheorem{theorem}{Theorem}[section]
\newtheorem{lemma}[theorem]{Lemma}
\newtheorem{proposition}[theorem]{Proposition}
\theoremstyle{definition}
\newtheorem{definition}[theorem]{Definition}
\theoremstyle{remark}
\newtheorem{remark}[theorem]{Remark}

\date{}

\newcommand{\manuscripttitle}[1]{\title[#1AP-SLDG schemes for Boltzmann]{#1Asymptotic-preserving semi-Lagrangian discontinuous Galerkin schemes for the Boltzmann equation}}
\author[X. Cai]{Xiaofeng Cai}
\address{Research Center of Mathematics, Faculty of Arts and Sciences, Beijing Normal University, Zhuhai, 519087, P.R. China, and Guangdong-Hong Kong-Macao Joint Laboratory for Mathematical Technology Research and Development,Beijing Normal University, Zhuhai 519087, China}
\email{xfcai@bnu.edu.cn}
\author[Z. Hao]{Zhen Hao}
\address{School of Mathematics and Statistics, Wuhan University, Wuhan 430072, P.R. China}
\email{zhhao\_math@whu.edu.cn}
\author[L. Liu]{Liu Liu}
\address{The Chinese University of Hong Kong, Hong Kong}
\email{lliu@math.cuhk.edu.hk}
\author[J. Wan]{Jiayu Wan}
\address{ShanghaiTech University, Shanghai 201210, China}
\email{wanjy@shanghaitech.edu.cn}
\thanks{The work of the third author was supported by National Key R\&D Program of China (2021YFA1001200), Ministry of Science and Technology in China, General Research Fund (14303022,  14301423, 14307125 \& 14310726) funded by Research Grants Council of Hong Kong.
The work of the first author was partially supported by the NSFC (No 12671463 and 12201052), Guangdong basic and applied basic research foundation[2025A1515012182], Guangdong and Hong Kong Universities “1+1+1” Joint Research Collaboration Scheme [2025A0505000014],Guangdong Major Project of Basic and Applied Basic Research (2021B0301030007).}

\manuscripttitle{}

\begin{document}
\begin{abstract}
We present an asymptotic-preserving semi-Lagrangian discontinuous Galerkin scheme for the Boltzmann equation that effectively handles multi-scale transport phenomena. The main challenge lies in finding the correct limit of the scheme as Knudsen number $\varepsilon \to 0$ in the semi-Lagrangian framework. We demonstrate how to overcome this difficulty in the paper, and our theoretical analysis establishes accuracy order conditions for both the IMEX-RK time integration and the limiting scheme. We also employ hypocoercivity techniques to establish stability for the linearized model. Numerical experiments for various test problems validate our proposed scheme's accuracy, asymptotic-preserving property, and robustness in various regimes, which demonstrates its effectiveness for multi-scale kinetic simulations.
\end{abstract}
\maketitle

\section{Introduction}

The Boltzmann equation provides a fundamental kinetic description of rarefied gas dynamics by modeling the evolution of the particle distribution function in phase space~\cite{Cercignani-Illner-Pulvirenti}. This equation is particularly important for simulating gas flows where continuum models such as the Navier-Stokes equations fail to capture the correct physics. The main computational difficulties of the Boltzmann equation often arise from the high-dimensionality of the phase space, and the nonlinear, nonlocal nature of the collision operator. The approaches for numerical solution of the Boltzmann equation date back to as early as 1940s~\cite{Grad-1949}. The most famous examples are the Direct Simulation Monte-Carlo (DSMC) methods by Bird~\cite{Bird} and by Nanbu~\cite{Nanbu-1980}. However, the stochastic nature of the DSMC
method introduces high statistical noise in low-speed flows, resulting in low accuracy. Over the years, deterministic methods, being free of statistical noise, have proven particularly advantageous for solving the Boltzmann equation, as reviewed in~\cite{Dimarco-Pareschi}. Among them, the Fourier spectral methods have been popularly used since the pioneering work~\cite{Pareschi-Russo}. Fast spectral methods were developed afterwards that can solve the collision operator in $\mathcal{O}(N^{d} \log N)$, see, for example, \cite{Mouhot-Pareschi,Filbet-Mouhot-Pareschi, Jaiswal-Alexeenko-Hu, Gamba-Haack-Hauck-Hu}.

Another fundamental component of the Boltzmann equation is the transport dynamics. Traditionally, transport dynamics have been treated using Eulerian finite difference and finite volume methods, such as upwind schemes~\cite{Godunov-1959} and weighted essentially non-oscillatory (WENO) schemes~\cite{Liu-Osher-Chan-1994}. Semi-Lagrangian (SL) schemes have gained popularity as an efficient method to treat the transport term across diverse applications including plasma simulations, climate modeling, and kinetic transport~\cite{Sonnendrucker-Roche-Bertrand-Ghizzo-1999, Lin-Rood-1996}. 
Different to Eulerian methods, which are based on a fixed spatial grid perspective and typically introduce Courant–Friedrichs–Lewy (CFL) time step restrictions for stability, SL methods avoid these restrictions by tracing back particle motions and integrating along characteristic trajectories, enabling larger time steps~\cite{Filbet-Sonnendrucker-2003}. Recently, \cite{Boscarino-Cho-Russo-2024} considered semi-Lagrangian methods applied to the Boltzmann equation. On the other hand, the discontinuous Galerkin (DG) method has been widely used in spatial discretization~\cite{Reed-Hill-1973}. Due to its conservation properties, $h$-$p$ adaptivity, and flexibility in handling complicated geometries, it has been preferred in the study of many stationary and non-stationary problems~\cite{Shu-2009, Cangiani-Georgoulis-Houston-2014}. In particular, we mention the nodal-DG framework~\cite{Hesthaven2007nodal, Crouseilles-Mehrenberger-Vecil-2011}, which offers flexibility in complex geometries and convenience when dealing with non-polynomial functions~\cite{Jaiswal-2022}. Combining the advantages of both the SL and DG methods has been a strategy pursued by many. For example, high-order conservative SLDG methods for kinetic and fluid equations have been developed in~\cite{Cai-Guo-Qiu-2017, Cai-Guo-Qiu-2019}. Recently, SLDG methods have 
also shown their advantages for simulation on GPUs~\cite{Einkemmer-Moriggl-2024}. 

In multi-scale scenarios, the stiffness of the collision operator presents additional computational bottlenecks. Severe time step restrictions like $\Delta t = O(\varepsilon)$ are imposed on explicit numerical methods. To overcome this challenge, asymptotic-preserving (AP) schemes have been a popular and robust computational approach for kinetic and hyperbolic problems~\cite{Jin-1999}. These schemes preserve the asymptotic transition from one scale to another at the discrete level, ensuring accuracy across regimes. The main advantage of AP schemes lies in their efficiency across all regimes, especially in the hydrodynamic or diffusive regimes. Unlike traditional methods, AP schemes do not need the numerical resolution of small physical parameters yet can still accurately capture the macroscopic behavior. In recent decades, AP schemes have found applications in various problems including semiconductor devices~\cite{Jin-Liu-2017}, radiative transfer equations~\cite{Fu-Cheng-Li-Xiong-Wang-2025}, and many others. For a comprehensive review, see \cite{AP2022}.
For the Boltzmann equation, a popular and effective approach is the BGK penalization method~\cite{Filbet-Jin}. For high-order temporal discretization, asymptotically accurate (AA) schemes that are able to preserve the order of accuracy in the fluid regime are discussed in~\cite{Dimarco-Pareschi-2013}.

Developing high-order SLDG schemes that satisfy the AP or AA property brings significant numerical challenges and often requires additional efforts. For instance, Ding et al. \cite{Ding-Qiu-Shu-2023} observed a reduction in order of numerical schemes when designing high-order AP-SLDG method for the BGK equation in the limiting regime. In particular, to ensure a third-order accuracy in this regime, additional condition must be imposed on the scheme. The Boltzmann equation, on the other hand, leads to greater challenges and difficulties. First, unlike in purely Eulerian methods, direct integration over velocity space fails to eliminate the stiff terms, rendering the derivation of the limiting scheme as $\varepsilon \to 0$ particularly challenging. Second, the temporal truncation error analysis is different from standard Eulerian method due to the coupling between temporal and spatial variables in the characteristic tracing operator. Therefore, the implicit-explicit Runge-Kutta (IMEX-RK) methods developed in the Eulerian framework \cite{Dimarco-Pareschi-2013} require extra attention when applied directly to the SL framework. Furthermore, other properties such as preserving positivity and ensuring numerical stability in the SL framework remain challenging and should be addressed.  

\textbf{Main contributions.}
In this work, we develop an AP-SLDG scheme for the Boltzmann equation with the derivation of the limiting scheme as the Knudsen number $\varepsilon \to 0$, and carefully investigate our semi-discretized scheme in aspects of order conditions, asymptotic accuracy, positivity preservation, and stability. To be specific, to obtain the correct limiting scheme, we introduce a microscopic auxiliary variable $g$ and derive a coupled evolution of the macroscopic states and the microscopic correction. Then, we study a class of IMEX-RK schemes in terms of their local truncation errors and asymptotic accuracy. One main difficulty in the error analysis under our SL setting arises from the involvement of time variable in the characteristic tracing operator, requiring a more careful Taylor expansion analysis. Moreover, we derive some sufficient conditions to maintain positivity of the distribution solution in the SL framework. We also adopt hypocoercivity tools to establish the stability of semi-discretized scheme for the linearized model, whose extension to the nonlinear Boltzmann case is deferred to a future work. Extensive numerical experiments are conducted to demonstrate the accuracy, AP property, and robustness of our proposed schemes in several test problems under different regimes. 

The rest of the paper is organized as follows. Section~\ref{sec:model} gives some background on the Boltzmann equation. Section~\ref{sec:SLDG} presents the 
AP-SLDG method and high-order IMEX-RK schemes. Section~\ref{derive limit scheme} is devoted to the derivation of the limiting scheme. Section~\ref{sec:analysis} provides theoretical analysis on positivity, order constraints and stability. In section ~\ref{sec:numerics}, we conduct various numerical tests to validate the accuracy and efficiency of our schemes. Section~\ref{sec:conclusion} shows the conclusion and future work.

\section{The Boltzmann equation}
\label{sec:model}
We consider the Boltzmann equation
\vspace{-4pt}
{\small
\begin{equation}\label{eq:BoltzmannEquation}
    \partial_t f + v \cdot \nabla_{x} f = \frac{1}{\mathrm{\varepsilon}} Q(f), 
\end{equation}
}where $f(x,v,t)$ is the probability density function of particles at position $x\in\Omega_{x} \subset \mathbb{R}^{d_x}$ with velocity $v\in\mathbb{R}^{d_v}$ at time $t>0$. 
The Knudsen number $\varepsilon$ is the dimensionless mean free path, defined as the average distance
between two collisions and characterizes the degree of rarefiedness of the gas. The operator $Q(f)$ describes particle collisions and act on the velocity variable. In particular, 
{\small
\begin{equation}
\label{eq:BoltzmannOperator}
    Q(f)(v) = \int_{\mathbb{R}^{d_v}} \int_{\mathbb{S}^{d_v-1}} \mathcal{B}\left(\left|v - v_*\right|, \sigma\right) \left[f(v') f(v_*') - f(v) f(v_*)\right] \, \mathrm{d}\sigma \, \mathrm{d}v_*,
\end{equation}
}where $(v, v_*)$ and $(v', v_*')$ are the pre-collision and post-collision velocity pairs related by conservation of momentum and energy
{\small
\begin{equation}\label{Boltz_eqn}
    v' = \frac{v + v_*}{2} + \frac{\left|v - v_*\right|}{2} \sigma, \quad
    v_*' = \frac{v + v_*}{2} - \frac{\left|v - v_*\right|}{2} \sigma,
\end{equation}
}with $\sigma$ varying over $\mathbb{S}^{d_v-1}$.
The quantity $B(\geq 0)$ is the collision kernel depending
only on $|v - v_*|$ and the scattering angle $\sigma$. 
In this work, we concern ourselves with the so-called Maxwell molecules with $d_v=2$~\cite{Villani_review}.

Throughout the paper, we use the notation $\langle \cdot,\cdot \rangle$ 
to stand for the standard $L^2$ inner product on the velocity space $L^2(\mathbb{R}^{d_v})$. The density, bulk velocity, and temperature of the distribution function $f$ are defined as
{\small
    $\rho = \langle f, 1 \rangle, \; u = \frac{1}{\rho} \langle f, v \rangle, \; T = \frac{1}{3\rho} \langle f, |v - u|^2 \rangle\,.$ 
}We use $U$ to represent the \textit{moments}, composed of the mass, momentum and energy. The moments of $f$ are
{\small
\begin{equation}
\label{eq:MacroObservables}
U := \left( \rho, \rho u, E \right)^\intercal = \left(\rho, \; \rho u, \; 3 \rho T + \rho |u|^2\right)^\intercal = \langle f, \phi \rangle,
\end{equation} 
}where $\phi(v)=(1,v,|v|^2)^\intercal$. The collision operator $Q$ satisfies the mass, momentum, and energy conservation
$\langle Q, \phi \rangle = 0.$
Hence $\phi$ are called the collision invariants~\cite{Cercignani}.
The model \eqref{eq:BoltzmannEquation} reaches the thermodynamic equilibrium when $f$ is a Maxwellian distribution~\cite{Cercignani} given by 
{\small
\begin{equation}
    \label{eq:Maxwellian}
     \mathcal{M}_f = \frac{\rho}{(2 \pi T)^{d_v/2}} 
     \exp\left(-\frac{|v - u|^2}{2T}\right).
\end{equation}
}Here $\mathcal{M}_f$ is called \textit{the Maxwellian associated with $f$}. Later we will use the notation $\mathcal{M}[U]$ to denote the Maxwellian constructed from the moments $U$.

\vspace{5pt}

\textbf{Fluid limit.} Integrating~\eqref{eq:BoltzmannEquation} against the collision invariants in the velocity space leads to the unclosed conservations laws
{\small
\begin{equation*}
    \partial_t \langle f \phi \rangle + \nabla_{x} \cdot \langle v \phi f \rangle = 0.
\end{equation*}
}While as $\varepsilon \rightarrow 0$, $f$ will be relaxed to the equilibrium state $\mathcal{M}_f$~\cite{Bardos-Golse-Levermore}. Thus,
we have recovered the limiting Euler equation
{\small
\begin{equation}
\label{eq:EulerEquation}
    \partial_t U = \mathcal{T} (U), 
\end{equation}
}where {\small $U = \langle \phi \mathcal{M}_f \rangle = (\rho, \rho u, E)^\intercal$}, 
{\small $\mathcal{T} (U) = -\nabla_x \cdot \langle v \phi \mathcal{M}_f \rangle = - (\rho u, \rho u \otimes u + \rho T I, E u + \rho T u)^\intercal$}.

\section{The nodal AP-SLDG method for the Boltzmann equation}

\label{sec:SLDG}

In this section, we will introduce the AP-SLDG method for the Boltzmann equation and will give the fully discretized scheme. 

\subsection{Nodal DG framework}
\label{sec:nodalDG}
We consider the nodal DG method in this work.
The advantages of dealing with nodal values as pointed out in~\cite{Ding-Qiu-Shu-2024}, are the flexibility to work with Maxwellian distributions~\eqref{eq:Maxwellian},
in addition to solving the problem with spatially dependent Knudsen numbers $\varepsilon(x)$, as shown in our last test. 
In particular, we adopt the nodal DG formulation developed in~\cite{Crouseilles-Mehrenberger-Vecil-2011, Einkemmer-Moriggl-2024}. 

We consider an one-dimensional spatial domain $\Omega_x = [x_L, x_R]$ and partition it into $N_x$ elements
\(
    x_L = x_{\frac{1}{2}} < x_{\frac{3}{2}} < ... < x_{N_x + \frac{1}{2}} = x_R
\).
Each element is denoted by $I_j = [x_{j - \frac{1}{2}}, x_{j + \frac{1}{2}}]$ with interval length $\Delta x_j =  x_{j + \frac{1}{2}} -  x_{j - \frac{1}{2}}$.
Numerical solutions and test functions are considered in the approximation space
\(
    V^k = \{ p : p|_{I_j} \in \mathcal{P}^k(I_j),\ j = 1, ..., N_x \},  
\)
where $k\in\mathbb{N}$ and $\mathcal{P}^k(I_j)$ denotes the set of polynomials of degree at most $k$ over the interval $I_j$. 
In each interval, we choose $k + 1$ Gauss-Legendre points as nodal points. The Lagrangian polynomials that interpolate at these nodal points are used as basis functions. The numerical solution in this work is represented by
{\small
\begin{equation}
    \label{eq:solution_representation}
    f(x,v) = \sum_{j = 1}^{N_x} \sum_{p_j = 0}^{k} f(x_{j, p_j}, v) \ell_{j, p_j}(x),
\end{equation}
}where $x_{j, p_j}$ is the $p_j$-th Gauss-Legendre point scaled to cell $j$, $\ell_{j, p_j}$ is the Lagrange basis function corresponding to $x_{j, p_j}$, and $k$ is the highest degree of basis polynomials.

\subsection{Characteristic Galerkin weak formulation}
We employ the weak formulation of characteristic Galerkin method proposed in~\cite{Celia-Russel-Herrera-Ewing-1990}. Consider the Boltzmann equation~\eqref{eq:BoltzmannEquation} in the 1D spatial variable case, with $x\in\Omega_x$. Let $t\in[0, T]$, the characteristic form of~\eqref{eq:BoltzmannEquation} can be written as
{\small
\begin{equation}
    \label{eq:BoltzmannEquationODE}
            \frac{\mathrm{d}}{\mathrm{d}t} f(t, x(t)) = \frac{1}{\varepsilon} Q(f), \quad 
            \frac{\mathrm{d}}{\mathrm{d}t} x(t) = v.
\end{equation}
}To simplify notations, here $v$ denotes the first component of the velocity variable. For each $v$, the characteristic curve is a straight line with slope $v$. For any $\Psi\in L^2(\Omega_x)$, we consider the following adjoint final-value problem as in~\cite{Cai-Guo-Qiu-2017}
{\small
\begin{equation}
    \label{eq:adjoint_problem}
            \partial_t \psi(x,t) + \partial_x (v \psi(x,t)) = 0, \quad 
            \psi(t=T) = \Psi(x),  
\end{equation}
}
the solution to which is simply $\psi(t,x) = \Psi(x + v (T-t))$ for $t\in[0,T]$.

\begin{proposition}
    \label{prop:weak_formulation}
Let $I$ be a line segment on the interval $\Omega_x$, assume $f$ and $\psi$ are solutions of~\eqref{eq:BoltzmannEquationODE} and~\eqref{eq:adjoint_problem} respectively, then the following identity holds
{\small\begin{equation}
    \label{eq:weak_formulation}
    \frac{\mathrm{d}}{\mathrm{d}t} \int_{I(t)} f \psi \mathrm{d}x =  \int_{I(t)} \frac{1}{\varepsilon} Q(f) \psi \mathrm{d}x.
\end{equation}}
\begin{proof}
    This identity can be proved by applying Leibniz rule and divergence theorem and using~\eqref{eq:BoltzmannEquation} and~\eqref{eq:adjoint_problem}.
\end{proof}
\end{proposition}
Proposition~\ref{prop:weak_formulation} leads to the following characteristic Galerkin weak formulation of the Boltzmann equation~\eqref{eq:BoltzmannEquationODE}. 
\begin{proposition}
Given any time interval $[t_1, t_2]\subset [0, T]$, line segment $[x_l, x_r] \subset \Omega_x$, and velocity $v$. Let $f$ be a solution to the Boltzmann equation \eqref{eq:BoltzmannEquationODE} with initial condition $f(t_1, x, v) = f_0(x, v)$, and $\psi$ be a solution to the adjoint problem \eqref{eq:adjoint_problem} with final condition $\psi(t_2, x) = \Psi(x)$, then we have
    {\small
\begin{equation}
    \label{eq:characteristic_tracing}
    \begin{aligned}
        \int_{x_l}^{x_r} f(t_2, x, v) & \Psi(x) \mathrm{d}x 
        = \int_{x_l - v \Delta t}^{x_r - v \Delta t} f_0(x, v) \Psi(x + v \Delta t) \mathrm{d}x \\
        &+ \int_{t_1}^{t_2} \int_{x_l - v (t_2-t)}^{x_r - v (t_2-t)} \frac{1}{\varepsilon} Q(f(t, x, v)) \Psi(x + v (t_2-t)) \mathrm{d}x\mathrm{d}t,
    \end{aligned}
\end{equation}
}
with $\Delta t = t_2 - t_1$.
\end{proposition}

\subsection{The AP-SLDG scheme}

We develop the AP-SLDG discretization by the following steps: (1) time discretization and BGK penalization based on~\eqref{eq:characteristic_tracing}; (2) spatial discretization by applying the nodal DG representation in~\eqref{eq:solution_representation}.

\subsubsection*{Step 1: Temporal discretization and BGK penalization}

The discretization of the Boltzmann equation~\eqref{eq:BoltzmannEquation} in the SLDG framework will largely be based on~\eqref{eq:characteristic_tracing}. 
In this part, an asymptotic-preserving time discretization is studied. 
For notation convenience, we denote $f^n(x, v) = f(t^n, x, v)$ as the numerical solution at time $t^n$, $U^{n} = U_{f^n}$, and $\mathcal{M}^n=\mathcal{M}_{f^n}$ according to~\eqref{eq:MacroObservables} and~\eqref{eq:Maxwellian}. 
To deal with the stiffness in the last integral, the collision operator needs an implicit treatment yet can be solved efficiently, a core idea in designing asymptotic-preserving schemes for the Boltzmann equation~\cite{Dimarco-Pareschi, Filbet-Jin}.
We adopt the BGK penalization approach proposed in~\cite{Filbet-Jin} that features the following decomposition
$Q(f) = Q(f) - \beta Q_{BGK}(f) + \beta Q_{BGK}(f).$
where $Q_{BGK}(f) = \mathcal{M}_f - f$
denotes the BGK operator, and $\mathcal{M}_f$ is the Maxwellian associated with $f$. The stiff part is treated implicitly, while the less stiff part is treated explicitly. 
The penalization parameter is chosen as 
$\beta = \max_v Q^{-}(f) := \max_v \int_{\mathbb{R}^{d_v}} \int_{\mathbb{S}^{d_v-1}} \mathcal{B}\left(\left|v - v_*\right|, \sigma\right)  f(v_*) \, \mathrm{d}\sigma \, \mathrm{d}v_*.$
Let $[t_1, t_2] = [t^n, t^{n+1}]$, $f_0=f^n$, and $\Psi(x)\in V^k$ in~\eqref{eq:characteristic_tracing}. This leads to the following one-step time discretization
{\small
\begin{equation}
\label{eq:one_step_discretization}
\begin{aligned}
    \int_{x_{j-1/2}}^{x_{j+1/2}} & f^{n+1}(x, v)  \Psi(x) \mathrm{d}x = \int_{x_{j-1/2} - v \Delta t}^{x_{j+1/2} - v \Delta t} f^{n}(x, v) \Psi(x + v \Delta t) \mathrm{d}x 
      \\ &  + \Delta t \int_{x_{j-1/2} - v \Delta t}^{x_{j+1/2} - v \Delta t} \frac{1}{\varepsilon} \left( Q(f^n) - \beta^n Q_{BGK}(f^n) \right)(x, v) \Psi(x + v \Delta t) \mathrm{d}x \\ & + \Delta t \int_{x_{j-1/2}}^{x_{j+1/2}} \frac{1}{\varepsilon} \beta^{n+1} Q_{BGK}(f^{n+1})(x, v) \Psi(x) \mathrm{d}x,
\end{aligned}
\end{equation}
}where $I(t^{n+1})\subset \Omega_x$ is some interval, and $I(t^n) = I(t^{n+1}) - v \Delta t$ is its upstream interval at time $t^n$.

\subsubsection*{Step 2: DG spatial discretization}
In this section, we introduce the spatial discretization of~\eqref{eq:one_step_discretization} by incorporating the DG representation in~\eqref{eq:solution_representation}. To simplify the discussion, we define the following auxiliary variables representing numerical values at different time levels
{\small
\begin{align*}
    \Theta^{n+1} = f^{n+1} - \Delta t \frac{\beta^{n+1}}{\varepsilon} (\mathcal{M}^{n+1} - f^{n+1}), \; 
    \Xi^{n} = f^n + \frac{\Delta t}{\varepsilon} (Q(f^n) - \beta^n (\mathcal{M}^n - f^n)). 
\end{align*}
}Using the Lagrange polynomials $\Psi(x) = \ell_{j, p_j}(x)$ as test functions and inserting the DG representation~\eqref{eq:solution_representation} into~\eqref{eq:one_step_discretization}, we perform Gauss quadrature to evaluate the integrals. Through a series of change of variables and quadrature evaluations (see appendix A), this leads to a matrix multiplication form
{\small
\begin{align}
\label{eq:matrix_form}
\Theta^{n+1}_{j, \cdot} &= A(\alpha) \Xi^n_{j^\ast, \cdot} + B(\alpha) \Xi^n_{j^\ast + 1, \cdot},
\end{align}
}where the matrices $A(\alpha)$ and $B(\alpha)$ are defined as
{\small
\begin{align*}
    A(\alpha)_{p_j, p} &= \omega_{p_j}^{-1} (1 - \alpha) \sum_q \omega_q \tilde{\ell}_p (\alpha + u_q (1 - \alpha)) \tilde{\ell}_{p_j} (u_q (1 - \alpha)), \\
    B(\alpha)_{p_j, p} &= \omega_{p_j}^{-1} \alpha \sum_q \omega_q \tilde{\ell}_p (\alpha u_q) \tilde{\ell}_{p_j} (\alpha(u_q - 1) + 1)),
\end{align*}
}with $\tilde{\ell}_{p}(s)$ representing Lagrange polynomials on $[0,1]$, $\omega_q$ and $u_q$ the Gauss-Legendre weights and nodes, and $j^\ast$ and $\alpha$ determined by the upstream tracing condition $x_{j - 1/2} - v \Delta t = x_{j^{\ast} - 1/2} + \alpha \Delta x$, $\alpha\in [0,1)$. 

Converting back to the original variables using the definitions of $\Theta^{n+1}$ and $\Xi^n$, we can express the scheme in terms of a reconstruction operator. For a given numerical function $f$ with nodal values $\{f_{j, \cdot}\}_{j = 1}^{N_x}$, we define $\mathcal{S}_{v, \Delta t}[f]$ as the numerical function with nodal values
{\small $\mathcal{S}_{v, \Delta t}[f]_{j, \cdot} := A(\alpha_v) f_{j^\ast_v, \cdot} + B(\alpha_v) f_{j^\ast_v + 1, \cdot}, \quad j = 1, \cdots, N_x.$
} The Eulerian AP-SLDG scheme is then given by
{\small
\begin{equation}
    \label{eq:AP-SLDG}
    \begin{aligned}
        f^{n+1} &= \mathcal{S}_{v, \Delta t}[f^n] + \mathcal{S}_{v, \Delta t}[\frac{\Delta t}{\varepsilon} (Q(f^n) - \beta^n (\mathcal{M}^n - f^n))] + \Delta t \frac{\beta^{n+1}}{\varepsilon} (\mathcal{M}^{n+1} - f^{n+1}).
    \end{aligned}
\end{equation}}
In fact, $\mathcal{S}_{v, \Delta t}[f](x, v)$ is a reconstruction of $f(x-v\Delta t, v)$.

\subsection{High-order AP-SLDG schemes}
\label{sec:IMEXRK}

In the DG framework, the spatial order can be improved by adjusting the polynomial order $k$. To achieve high-order temporal accuracy, we need to construct IMEX-RK schemes. To simplify notations, we define $G_{P}(f)=Q(f) - \beta(f) Q_{BGK}(f), \; Q_{P}(f)=\beta(f) Q_{BGK}(f).$
We consider the $s$-stage IMEX-RK method characterized by the following double Butcher tableau~\cite{Ascher-Ruuth-Spiteri-1997}
\begin{equation}
\label{eq:ButcherTable}
\begin{aligned}
\begin{array}{cc}
    \begin{array}{c|c}
        c \vphantom{\tilde{c}} & A \vphantom{\tilde{A}} \\
        \hline
        & \rule{0pt}{15pt} b^T \vphantom{\tilde{b}}
    \end{array} & \hspace{15pt}
    \begin{array}{c|c}
        \tilde{c} & \tilde{A} \\
        \hline
        & \rule{0pt}{15pt} \tilde{b}^T,
    \end{array}
\end{array}
\end{aligned}
\end{equation}
where $A = (a_{ij})$ and $\tilde{A} = (\tilde{a}_{ij})$ are $s\times s$ lower triangular matrices with $a_{ij} = 0$ for $j \geq i$, and $\tilde{a}_{ij} = 0$ for $j > i$. The coefficients $c_i$ and $\tilde{c}_i$ are given by
$c_i = \sum_{j=1}^{i} a_{ij}, \quad \tilde{c}_i = \sum_{j=1}^{i} \tilde{a}_{ij}. $
Here $b = (b_i)$ and $\tilde{b} = (\tilde{b}_i)$ represent the quadrature weights for internal stages of the RK method. In~\eqref{eq:ButcherTable}, the table on the left corresponds to the `explicit part' of the scheme, while that on the right corresponds to the `implicit part'.

Since the problem is stiff, we always consider globally-stiffly-accurate (GSA) methods~\cite{Hairer-Wanner-1996}, see Definition~\ref{def GSA}. In particular, GSA methods mean that the next step solution is always equal to the last stage solution, namely 
$f^{n+1} = f^{(s)}$.
Since there are several stages involved, we make the simplification of notations
$\mathcal{S}_{i,j} := \mathcal{S}_{v,(c_{i}-c_{j})\Delta t}$, $\tilde{\mathcal{S}}_{i,j} := \mathcal{S}_{v,(\tilde{c}_{i}-\tilde{c}_{j})\Delta t}$.  
Using the Butcher tableaux in~\eqref{eq:ButcherTable} to treat the implicit part and the explicit part of~\eqref{eq:AP-SLDG}, 
we derive the IMEX-RK method for the AP-SLDG scheme
{\small
\begin{equation}\label{eq:IMEX-RK}
\begin{aligned}
    f^{(i)} &= \tilde{\mathcal{S}}_{i,0} [f^n] + \Delta t \sum_{j=1}^{i-1} a_{ij} \mathcal{S}_{i,j} [\frac{1}{\varepsilon}G_{P}(f^{(j)})] + \Delta t\sum_{j=1}^{i} \tilde{a}_{ij} \tilde{\mathcal{S}}_{i,j} [\frac{1}{\varepsilon}Q_{P}(f^{(j)})], \\
    \end{aligned}
\end{equation}
}where $i\in\{1,\cdots,s\}$ denotes the $i$-th internal stage. 

We list out a few Butcher tableaux~\eqref{eq:IMEX-RK-1st}, \eqref{eq:IMEX-RK-DP2A242}, \eqref{eq:IMEX-RK-ARS443} in Appendix D. These methods will be denoted as \texttt{FBEuler}, \texttt{DP2A242}, and \texttt{ARS443}, respectively.
Note that the method \texttt{FBEuler} corresponds to the scheme~\eqref{eq:AP-SLDG}.

\begin{remark}
Unlike in Eulerian frameworks, in~\eqref{eq:IMEX-RK}, the stage abscissae have to be specified in the shifting of $f^n$ for each stage since $\tilde{c} \neq c$ for many Butcher tableau. However, for our purpose, we will mainly work on CK schemes(see Appendix B for definition), where the coefficients can be chosen so that $\tilde{c}=c$. Therefore, it does not matter for us if $\tilde{\mathcal{S}}_{i,0} [f^n]$ or $\mathcal{S}_{i,0} [f^n]$ is chosen in \eqref{eq:IMEX-RK}. 
\end{remark}

\eqref{eq:IMEX-RK}, together with the DG spatial reconstruction discussed in section 3.3, forms the main part of our AP-SLDG scheme and constitutes the core of this paper.

\section{Derivation of the limiting scheme}\label{derive limit scheme}
In the study of asymptotic preserving schemes of kinetic equations, it is important to derive the limit of our scheme as $\varepsilon \to 0$ and analyze if the properties of the original scheme are preserved asymptotically(e.g. order of accuracy, stability, positivity etc.). In this section, we derive the limiting scheme of \eqref{eq:IMEX-RK} as $\varepsilon \to 0$, and we shall analysis the properties of the limiting scheme in section 5. 

In the Eulerian framework~\cite{Filbet-Jin}, a typical way to obtain the limiting scheme is to multiply the scheme by $\phi(v) = (1,v,|v|^2)^T$ and integrate in $v$. From the conservation properties of the collision operators, the stiff terms usually vanish and a limiting scheme of the limiting equation(Euler equations in our case) can be obtained by passing $\varepsilon \to 0$. However, in the SL framework, the situation is different. For example, if we apply the same strategy to the first-order scheme~\eqref{eq:AP-SLDG}, we would obtain
$$
U^{n+1} = \langle \mathcal{S}_{v,\Delta t}[f^n], \phi(v) \rangle + \Delta t \Big\langle \mathcal{S}_{v, \Delta t}[\frac{1}{\varepsilon} \left(Q(f^n) - \beta^n (M^n - f^n)\right)] , \phi(v) \Big\rangle\,,
$$
in which the stiff collision terms does not vanish because of the presence of $v$ in $\mathcal{S}_{v,\Delta t}$. This leads to additional challenges when deriving the limiting scheme in the SL framework. We shall resolve this issue in this section.

To begin with, we discuss the asymptotic-preserving(AP) property of the proposed scheme \eqref{eq:IMEX-RK}. We first review the following definition~\cite{Dimarco-Pareschi}. 
\begin{definition}
We call an IMEX-RK scheme strongly AP if \[ \lim_{\varepsilon \to 0} f^{n+1}  = \mathcal{M}[f^{n+1}] \] for any initial data $f^n$. We call the scheme weakly AP if the same equation holds when the initial data $f^n$ is well-prepared, ie: $f^n = \mathcal{M}[f^n] + O(\varepsilon)$. 
\end{definition}

The following theorem can be proved. 
\begin{theorem}\label{AP property for schemes}
Assume the IMEX-RK scheme \eqref{eq:IMEX-RK} is GSA. Then, if the scheme is of type A, it is strongly AP; if the scheme is of type CK, it is weakly AP.
\end{theorem}

\begin{proof}
We first assume the scheme is GSA and of type A.  From \eqref{eq:IMEX-RK}, we have 
{\small
\begin{align*}
    \varepsilon f^{(i)} - \varepsilon \mathcal{S}_{i,0} [f^n] = \Delta t \sum_{j=1}^{i-1} a_{ij} \mathcal{S}_{i,j} [G_{P}(f^{(j)})] + \Delta t \sum_{j=1}^{i} \tilde{a}_{ij} \tilde{\mathcal{S}}_{i,j} [Q_{P}(f^{(j)})]. 
\end{align*}
}Take $\lim_{\varepsilon \to 0}$ on both sides of the equation, we have 
{\small
\begin{equation}\label{limit equation for Q_P_fi}
    \tilde{a}_{ii} Q_{P}(f^{(i)}) = \sum_{j=1}^{i-1} a_{ij} \mathcal{S}_{i,j} [G_{P}(f^{(j)})] + \sum_{j=1}^{i-1} \tilde{a}_{ij} \tilde{\mathcal{S}}_{i,j} [Q_{P}(f^{(j)})]. 
\end{equation}
}Clearly as $\varepsilon \to 0$,  $\tilde{a}_{11} Q_{P}(f^{(1)})=0$ and hence $Q_{P}(f^{(1)})=0$ due to the invertibility of $\tilde{A}$. This shows $f^{(1)}=\mathcal{M}[f^{(1)}]$ as $\varepsilon \to 0$. 

Now we make an induction argument. Suppose $i \geq 2$ and $\lim_{\varepsilon \to 0} f^{(j)} = \mathcal{M}[f^{(j)}]$ for all $j < i$. Let $\varepsilon \to 0$ and consider \eqref{limit equation for Q_P_fi}. The right-hand-side of the equation vanishes due to the inductive hypothesis and the fact that local Maxwellians lie in the kernels of both $G_P$ and $Q_P$. Then $f^{(i)}=\mathcal{M}[f^{(i)}]$ for exactly the same reason as $f^{(1)}$. This completes the induction and shows that $f^{(i)}=\mathcal{M}[f^{(i)}]$, $\forall 1 \leq i \leq s$ as $\varepsilon \to 0$. In particular, since we assume the scheme is GSA, $f^{n+1} = f^{(s)} = \mathcal{M}[f^{(s)}]$ as $\varepsilon \to 0$ and hence the scheme is strongly AP. 

We assume the scheme is GSA and of type CK. Suppose that the initial data is well-prepared, then $f^{(1)} = f^{n} = \mathcal{M}[f^{n}]+O(\varepsilon)$. Hence \( \lim_{\varepsilon \to 0} f^{(1)} = \mathcal{M}[f^{(1)}]. \) Then using the fact that $\tilde{\hat{A}}$ is invertible, one can conduct the same induction argument as in the type A case to conclude that $f^{n+1} = f^{(s)} = \mathcal{M}[f^{(s)}]$ as $\varepsilon \to 0$, hence the scheme is weakly AP. 
\end{proof}

\begin{remark}\label{AP for intermediate stages}\raggedright
In fact, in the proof of the theorem above, we showed that \mbox{$\lim_{\varepsilon \to 0} f^{(i)}=\mathcal{M}[f^{(i)}]$} for all the intermediate stages $i$. 
\end{remark}

We now derive the limiting scheme of \eqref{eq:IMEX-RK} as $\varepsilon \to 0$. The key idea is that we introduce a microscopic auxiliary variable $g$ through the near-equilibrium decomposition $f=\mathcal{M}[U]+\varepsilon g$, with $\langle\phi g\rangle=0$. This yields a closed limiting scheme for the coupled evolution of the macroscopic variables $U$ and the microscopic correction $g$. Although the limiting scheme retains the microscopic auxiliary variable $g$, its macroscopic update is formally consistent with the Euler equations \eqref{eq:EulerEquation}, and a systematic analysis of its temporal order of accuracy will be presented subsequently in section \ref{sec:accuracy order analysis}. For clarity, throughout the rest of this subsection, we use the subscript $\varepsilon$ to denote the dependency of the variables and functions on $\varepsilon$. The time step $h=\Delta t>0$ is fixed as $\varepsilon \to 0$.

We shall restrict our attention to the situation where a well-prepared family of initial data is provided, i.e. 
\begin{equation}\label{assumption on well-prepared initial data}
f_\varepsilon^0=\mathcal{M}[U_\varepsilon^0]+\varepsilon g_\varepsilon^0,\qquad
\langle\phi g_\varepsilon^0\rangle=0, \qquad \sup_{\varepsilon}\| g_{\varepsilon}^{0}\|_{\mathcal{H}} \le C,\qquad \text{for } C \text{ independent of }\varepsilon,
\end{equation}
with
\begin{equation}\label{assumption on convergence of initial data}
U_\varepsilon^0\longrightarrow U^0,\qquad g_\varepsilon^0\longrightarrow g^0, \qquad  \| g^{0}\|_{\mathcal{H}} \le C, \qquad \text{for } \varepsilon \to 0, 
\end{equation}
where $\mathcal{H}=L^2(\mathbb T_x^d\times\mathbb R_v^d)$. Therefore, the limiting scheme shall be initialized with $(U^0,g^0)$. For every time stage under consideration, write
\begin{equation}\label{eq:discrete-decomposition}
\begin{aligned}
f_\varepsilon^n&=M[U_\varepsilon^n]+\varepsilon g_\varepsilon^n,
&U_\varepsilon^n&=\langle\phi f_\varepsilon^n\rangle,\\
f_\varepsilon^{(i)}&=M[U_\varepsilon^{(i)}]+\varepsilon g_\varepsilon^{(i)},
&U_\varepsilon^{(i)}
&=\langle\phi f_\varepsilon^{(i)}\rangle.
\end{aligned}
\end{equation}
These definitions imply $\langle\phi g_\varepsilon^n\rangle=\langle\phi g_\varepsilon^{(i)}\rangle=0$. For a given macroscopic state $U$, define
\begin{equation}\label{eq:local-linearization}
\mathcal L_U(g):=Q(\mathcal{M}[U],g)+Q(g,\mathcal{M}[U]).
\end{equation}
Notice that although $\mathcal L_U$ is linear in $g$, its coefficients
depend on the macroscopic state $U$. We further define the following quantity 
\begin{equation}\label{eq:finite-eps-R}
\begin{aligned}
R_\varepsilon^{(i)}:=&\mathcal{S}_{i,0}[M[U_\varepsilon^n]+\varepsilon g_\varepsilon^n]+h\sum_{j=1}^{i-1}a_{ij}\mathcal{S}_{i,j}\Bigl[(\mathcal L_{U_\varepsilon^{(j)}}+\beta I)g_\varepsilon^{(j)}+\varepsilon Q(g_\varepsilon^{(j)},g_\varepsilon^{(j)})\Bigr] \\
&-h\beta\sum_{j=1}^{i-1}
\widetilde a_{ij}\mathcal{S}_{i,j}[g_\varepsilon^{(j)}].
\end{aligned}
\end{equation}

We now state and prove the limiting scheme of \eqref{eq:IMEX-RK}. For simplicity, we assume \eqref{eq:IMEX-RK} is of type CK so that the implicit and explicit time stages can be chosen to be the same. We further require that \eqref{eq:IMEX-RK} is GSA so that the final update $f^{n+1}$ equals the last stage value $f^{(s)}$. 

\begin{theorem}
Assume that \eqref{assumption on well-prepared initial data},\eqref{assumption on convergence of initial data} holds for the initial data and the scheme \eqref{eq:IMEX-RK} is of type CK and GSA. Then the limiting scheme of \eqref{eq:IMEX-RK} as $\varepsilon \to 0$ is the following coupled macro-micro scheme:
\begin{equation}\label{eq:limiting-scheme}
\begin{aligned}
U^{(1)}&=U^n,&g^{(1)}&=g^n,\\
U^{(i)}&=\langle\phi R^{(i)}\rangle,&g^{(i)}&=\frac{R^{(i)}-M[U^{(i)}]}{h\beta \widetilde{a}_{ii}},\qquad i=2,\ldots,s,\\
U^{n+1}&=U^{(s)},&g^{n+1}&=g^{(s)}.
\end{aligned}
\end{equation}
where 
\begin{equation}\label{eq:limiting-R}
\begin{aligned}
R^{(i)}:=\mathcal{S}_{i,0}[\mathcal{M}[U^n]]+h\sum_{j=1}^{i-1}a_{ij}\mathcal{S}_{i,j}[(\mathcal L_{U^{(j)}}+\beta I)g^{(j)}]-h\beta\sum_{j=1}^{i-1}\widetilde a_{ij}\mathcal{S}_{i,j}[g^{(j)}].
\end{aligned}
\end{equation}
\end{theorem}

\begin{proof}
Consider $f_{\varepsilon}=\mathcal{M}[U_{\varepsilon}]+\varepsilon g_{\varepsilon}$. By bilinearity of $Q$ and $Q(\mathcal{M}[U_{\varepsilon}],\mathcal{M}[U_{\varepsilon}])=0$, we have the identity
\begin{equation}
Q(\mathcal{M}[U_{\varepsilon}]+\varepsilon g_{\varepsilon},\mathcal{M}[U_{\varepsilon}]+\varepsilon g_{\varepsilon})=\varepsilon\mathcal L_U(g_{\varepsilon})+\varepsilon^2Q(g_{\varepsilon},g_{\varepsilon}).
\end{equation}
Since $\langle\phi g_{\varepsilon}\rangle=0$, $Q_P(\mathcal{M}[U_{\varepsilon}]+\varepsilon g_{\varepsilon})=-\varepsilon\beta g_{\varepsilon}$. Consequently,
\begin{equation}\label{eq:collision-expansions}
\begin{aligned}
\frac1\varepsilon Q_P(\mathcal{M}[U_{\varepsilon}]+\varepsilon g_{\varepsilon})&=-\beta g_{\varepsilon},\\
\frac1\varepsilon G_P(\mathcal{M}[U_{\varepsilon}]+\varepsilon g_{\varepsilon})&=(\mathcal L_U+\beta I)g_{\varepsilon}+\varepsilon Q(g_{\varepsilon},g_{\varepsilon}).
\end{aligned}
\end{equation}
Using \eqref{eq:collision-expansions} and the fact that $S_{i,i}=I$, the stage equation \eqref{eq:IMEX-RK} becomes
\begin{equation}\label{eq:finite-eps-reconstruction}
\mathcal{M}[U_\varepsilon^{(i)}]+(\varepsilon+h\beta \widetilde{a}_{ii})g_\varepsilon^{(i)}=R_\varepsilon^{(i)}.
\end{equation}
Taking collision-invariant moments gives $U_\varepsilon^{(i)}=\langle\phi R_\varepsilon^{(i)}\rangle$. Hence the exact stage reconstruction is
\begin{equation}\label{eq:finite-eps-micro-macro}
\begin{aligned}
U_\varepsilon^{(i)}&=\langle\phi R_\varepsilon^{(i)}\rangle,\\
g_\varepsilon^{(i)}&=\frac{R_\varepsilon^{(i)}-\mathcal{M}[U_\varepsilon^{(i)}]}{\varepsilon+h\beta \widetilde{a}_{ii}}.
\end{aligned}
\end{equation}
Assumption \eqref{assumption on well-prepared initial data}, together with Theorem \ref{AP property for schemes} and remark \ref{AP for intermediate stages}, implies that 
\begin{equation}
\mathcal{S}_{i,0}[\varepsilon g_\varepsilon^n]\longrightarrow0,\qquad\varepsilon
Q(g_\varepsilon^{(j)},g_\varepsilon^{(j)})\longrightarrow 0.
\end{equation}
Therefore, formally we have $R^{(i)}=\lim_{\varepsilon \to 0}R_{\varepsilon}^{(i)}$. Passing this to the limit in \eqref{eq:finite-eps-micro-macro} yields the desired scheme \eqref{eq:limiting-scheme}. The last line of \eqref{eq:limiting-scheme} follows from GSA assumption. We comment that assumption \eqref{assumption on convergence of initial data} is needed so that we have an initialization $(U^{0},g^{0})$ compatible with \eqref{eq:IMEX-RK} to start the algorithm \eqref{eq:limiting-scheme}. 
\end{proof}

\section{Numerical analysis for the scheme}
\label{sec:analysis}

This section analyzes the proposed scheme's positivity-preserving, temporal order, and stability properties. We consider only constant $\varepsilon$ (spatially dependent $\varepsilon(x)$ complicates analysis within operator $\mathcal{S}$) and focus on semi-discrete schemes with temporal discretization, where $\mathcal{S}_{v,\Delta t}[f]=f(x-v\Delta t,v)$ and collision operators $Q$, $Q_{P}$ maintain local conservation. Fully discrete analysis is deferred to future work.

We organize the analysis as follows: Section 4.1 derives sufficient conditions for positivity preservation in the SL framework. Section 4.2 analyzes local truncation accuracy of scheme \eqref{eq:IMEX-RK} and the limiting scheme \eqref{eq:limiting-scheme}. Section 4.3 establishes stability of the IMEX Euler scheme near global equilibrium using hypocoercivity tools, where nonlinear operator $Q$ is approximated by its linearization.

\subsection{Positivity-preserving property} 
Physically, the solution $f$ to the Boltzmann equation \eqref{eq:BoltzmannEquation} represents the particle distribution function, which must satisfy $f \geq 0$. Numerically, the inclusion of Maxwellians within the penalization term would amplify instabilities through temperature within the exponential function if $f$ becomes negative. Therefore, maintaining the positivity of the numerical approximation $f^n$ across all time steps is essential for both physical accuracy and algorithmic stability.
In this subsection, we will discuss the constraints on the RK coefficients to ensure the positivity of our scheme. Similar to~\cite{Dimarco-Pareschi-2013}, the key to the positivity analysis is the decomposition of operators into the following forms
{\small 
\begin{equation}
\label{decomp of Q_P G_P}
Q_P(f) = \beta(\mathcal{M}[f] - f), \; G_P(f) = P(f) - \beta \mathcal{M}[f],
\end{equation}
}where $P(f) = Q(f) + \beta f$, which is positive for $f\geq0$ if we choose $\beta>0$ properly~\cite{Filbet-Jin}. For simplicity, we consider IMEX-RK schemes of type CK that are GSA such that $c_i=\widetilde{c}_i$ for all $i$. Then the constraints on the coefficients to ensure the positivity of our scheme \eqref{eq:IMEX-RK} can be summarized as follows. 

\begin{theorem}\label{positivity constraints for RK}
Suppose that the scheme \eqref{eq:IMEX-RK} is of type CK and GSA such that $c_i=\widetilde{c}_i$ for all $i$. Assume $f^n \geq 0$. Then a sufficient condition to ensure $f^{n+1} \geq0$ is given by 
\begin{equation}\label{eq:ck-positivity-conditions}
\begin{aligned}
&\widetilde a_{ii}\ge0, ~~(K\mathbf e)_i\ge0,
&& i=1,\ldots,s,\\
&(KA)_{ij}\ge0,~~\bigl(K(\widetilde A-A)\bigr)_{ij}\ge0,
&& 1\le j<i\le s.
\end{aligned}
\end{equation}
where $K=(I+z\widetilde A)^{-1},~\mathbf e=(1,\ldots,1)^T$ and $z=\frac{\beta\Delta t}{\varepsilon}>0$.
\end{theorem}

\begin{proof}
Since $c_i=\widetilde{c}_i$, we have $\mathcal{S}_{i,j}=\widetilde{\mathcal{S}}_{i,j}$ for all $j \le i$. Then from \eqref{eq:IMEX-RK}, the stage equations are
$$
f^{(i)}=\mathcal{S}_{i,0}[f^n]+\frac{\Delta t}{\varepsilon} \sum_{j<i}a_{ij}\mathcal{S}_{i,j}[G_P(f^{(j)})]+\frac{\Delta t}{\varepsilon}\sum_{j\le i}\widetilde a_{ij}\mathcal{S}_{i,j}[Q_P(f^{(j)})].
$$
Recall the identities $Q_P(f) = \beta(\mathcal{M}[f] - f), ~G_P(f) = P(f) - \beta \mathcal{M}[f]$, where $f\geq 0$ implies $P(f)\geq 0$ if $\beta$ is properly chosen. Writing $M^{(i)}=M[f^{(i)}]$, $P^{(i)}=P(f^{(i)})$ and $z=\frac{\beta\Delta t}{\varepsilon}$, then the stage equations above can be rewritten as 
\begin{equation}\label{stage equation 2}
f^{(i)}+z\sum_{j\le i}\widetilde a_{ij}\mathcal{S}_{i,j}[f^{(j)}]=
\mathcal{S}_{i,0}[f^n]+\frac{z}{\beta}\sum_{j<i}a_{ij}\mathcal{S}_{i,j}[P^{(j)}]+z\sum_{j\le i}(\widetilde a_{ij}-a_{ij})\mathcal{S}_{i,j}[M^{(j)}].
\end{equation}
Let $K=(I+z\widetilde A)^{-1},~\mathbf e=(1,\ldots,1)^T.$ If $\widetilde a_{ii}\ge0$, then $d_i>0$, and $K$ exists even though the CK matrix $\widetilde A$ itself is singular. Multiply row $k$ of \eqref{stage equation 2} by $K_{ik}\mathcal S_{i,k}$ and sum over $k\le i$. Using $\mathcal S_{i,k}\mathcal S_{k,j}=\mathcal S_{i,j}$ and $K(I+z\widetilde A)=I$, we obtain
\begin{equation}\label{stage equation 3}
f^{(i)}=(K\mathbf e)_i\mathcal{S}_{i,0}[f^n]+\frac{z}{\beta}
\sum_{j<i}(KA)_{ij}\mathcal{S}_{i,j}[P^{(j)}]+z\sum_{j\le i}\bigl(K(\widetilde A-A)\bigr)_{ij}\mathcal{S}_{i,j}[M^{(j)}].
\end{equation}
Now define
\begin{equation}\label{def of G^i}
\begin{aligned}
G^{(i)}=& \mathcal{S}_{i,0}[f^n]+\frac{\Delta t}{\varepsilon} \sum_{j<i}a_{ij}\mathcal{S}_{i,j}[G_P(f^{(j)})]+\frac{\Delta t}{\varepsilon}\sum_{j<i}\widetilde a_{ij}\mathcal{S}_{i,j}[Q_P(f^{(j)})]\\
=&\mathcal{S}_{i,0}[f^n]+\frac{z}{\beta}\sum_{j<i}a_{ij}\mathcal{S}_{i,j}[P^{(j)}]+z\sum_{j<i}(\widetilde a_{ij}-a_{ij})\mathcal{S}_{i,j}[M^{(j)}]-z\sum_{j<i}\widetilde a_{ij}\mathcal{S}_{i,j}[f^{(j)}].
\end{aligned}
\end{equation}
Then $f^{(i)}=G^{(i)}+z\widetilde a_{ii}(M^{(i)}-f^{(i)}).$ Let $d_i=1+z\widetilde a_{ii}$, thus $f^{(i)}=\frac{G^{(i)}+z\widetilde a_{ii}M^{(i)}}{d_i}.$ The purpose of introducing $G^{(i)}$ is to establish its positivity before deriving the positivity of $M^{(i)}$ and $f^{(i)}$, as we now explain. 

Since $K_{ii}=1/d_i$ and $a_{ii}=0$, $\bigl(K(\widetilde A-A)\bigr)_{ii}=\frac{\widetilde a_{ii}}{d_i}.$
Multiplying \eqref{stage equation 3} by $d_i$ and
subtracting $z\widetilde a_{ii}M^{(i)}$ gives

\begin{equation}\label{G^i another form}
G^{(i)}=d_i(K\mathbf e)_i\mathcal{S}_{i,0}[f^n]+\frac{z}{\beta}d_i\sum_{j<i}(KA)_{ij}\mathcal{S}_{i,j}[P^{(j)}]+zd_i\sum_{j<i}\bigl(K(\widetilde A-A)\bigr)_{ij}\mathcal{S}_{i,j}[M^{(j)}].
\end{equation}
Importantly, the right-hand side of \eqref{G^i another form} contains only quantities associated previous stages. Now assume that the positivity conditions given in the statement of the theorem hold, i.e. 
\begin{equation}
\begin{aligned}
&\widetilde a_{ii}\ge0, ~~(K\mathbf e)_i\ge0,
&& i=1,\ldots,s,\\
&(KA)_{ij}\ge0,~~\bigl(K(\widetilde A-A)\bigr)_{ij}\ge0,
&& 1\le j<i\le s.
\end{aligned}
\end{equation}

For a CK scheme, the first rows of both tableaux vanish. Therefore, $c_1=0$ and
$f^{(1)}=f^n \geq 0$. Suppose now that $f^{(j)}\ge0$ for all $j<i$. Then the density(0th moment) and temperature(2nd moment) associated to $f^{(j)}$ are both positive, and hence $M^{(j)}\ge0$. We also have $P^{(j)}\ge0$ and $\mathcal{S}_{i,j}$ preserves positivity. Then \eqref{eq:ck-positivity-conditions}implies, through \eqref{G^i another form}, that $G^{(i)}\ge0.$

Since $f^{(i)}=G^{(i)}+z\widetilde a_{ii}(M^{(i)}-f^{(i)})$ and the BGK operator preserves collision invariants, we have $\int\psi(v)f^{(i)}\,dv=\int\psi(v)G^{(i)}\,dv.$ for $\psi(v)=\left(1,v,\frac{|v|^2}{2}\right)$. Therefore, $G^{(i)}\ge 0$ implies $M^{(i)}\ge 0$. Since $f^{(i)}=\frac{G^{(i)}+z\widetilde a_{ii}M^{(i)}}{d_i}$ and $\widetilde a_{ii}\ge0$, this yields $f^{(i)}\ge 0$. Since the scheme is GSA, $f^{n+1}$ equals the last stage value, and hence $f^{n+1}\geq 0$. This completes the proof. 
\end{proof}

\begin{remark}
We comment that \eqref{eq:ck-positivity-conditions} imposes sufficient conditions for the IMEX-RK coefficients to ensure the positivity of the scheme, not necessary ones. In particular, the two schemes tested in Section 6, \texttt{DP2A242} and \texttt{ARS443}, will not strictly satisfy \eqref{eq:ck-positivity-conditions}. But they do not suffer from positivity issues in most situations(except for \texttt{ARS443} in the Sod shock problem). We leave it to future work to find some higher-order IMEX-RK Butcher tableau that strictly satisfies the positivity restrictions stated in \eqref{eq:ck-positivity-conditions}. 
\end{remark}

\subsection{Conditions for the accuracy order} \label{sec:accuracy order analysis} 

In this subsection, we analyze the local temporal accuracy of the scheme~\eqref{eq:IMEX-RK} and the limiting scheme \eqref{eq:limiting-scheme} by deriving proper order conditions on the RK coefficients. The analysis is different from the Eulerian case since now the shifting operator $\mathcal{S}_{v,\Delta t}$ also contains the time step $\Delta t$. However, The shifting operator actually won't affect the temporal accuracy order for scheme \eqref{eq:IMEX-RK} with positive $\varepsilon$, as proved by Theorem \ref{thm:sl-imex-order} below. This is a big advantage of SL schemes. For the limiting scheme \eqref{eq:limiting-scheme}, an order reduction will typically be observed, which is common among numerical schemes for Boltzmann-type equations(see \cite{Ding-Qiu-Shu-2024},\cite{Dimarco-Pareschi}). We shall prove that \eqref{eq:limiting-scheme} is at least first-order consistent with the compressible Euler equations, and we also analyze sufficient conditions for the IMEX-RK coefficients to ensure second-order consistency of \eqref{eq:limiting-scheme}. Similar analysis for higher order accuracy of the limiting scheme can be performed as well, and we leave it for future study. 

\subsubsection{Order conditions for SL-IMEX-RK schemes with fixed $\varepsilon>0$}

The main result of this subsection is Theorem \ref{thm:sl-imex-order}, which essentially states that the presence of the shifting operator $\mathcal{S}_{v,\Delta t}$ in \eqref{eq:IMEX-RK} won't affect the standard IMEX-RK order conditions. To begin with, we provide a short introduction to order conditions for typical IMEX-RK schemes. 

Consider the autonomous equation 
\begin{equation}\label{autonomous y(t)}
y'=F(y)+G(y),
\end{equation}
where $y:\mathbb{R}\longrightarrow \mathbb{R}$, with $F$ treated explicitly by $(A,b)$ and $G$ implicitly by $(\widetilde A,\widetilde b)$. An IMEX-RK method has classical order $p$ if, starting from $y^n=y(t_n)$, its one-step error is
$O(h^{p+1})$, where $h=\Delta t$. The order conditions(on the coefficients $(A,b)$ and $(\widetilde A,\widetilde b)$) follow by matching the Taylor expansions of the exact and numerical solutions through degree $p$. For more references on the derivation of these conditions, we refer to \cite{Butcher2016}.

When the target space of $y$ remains continuous, we regard $u(t)=y(t,\cdot)$
as an element of a Banach space $\mathcal{H}$(in our case $\mathcal{H}=L^2(\mathbb{T}^d_{x}\times \mathbb{R}^d_{v})$). As long as certain regularity conditions hold, the algebraic order conditions are unchanged: we can still perform the Taylor expansions with ordinary derivatives replaced by Fr\'echet derivatives, and the local error is measured in the norm of $\mathcal{H}$. 

Now, consider the non-autonomous equation 
\begin{equation}\label{non-autonomous y(t)}
u'=F(t,u)+G(t,u).
\end{equation}
The equation can be written as the autonomous system
\begin{equation}\label{non-autonomous system}
\binom{\theta}{u}'=\binom{1}{F(\theta,u)}+\binom{0}{G(\theta,u)},\qquad \theta(t_n)=t_n.
\end{equation}
Applying the same IMEX--RK coefficients gives time stages $\Theta_i=t_n+h(A\mathbf 1)_i$. Under the common-node condition
$A\mathbf 1=\widetilde A\mathbf 1=c$, these are precisely the
prescribed evaluation times $t_n+c_i h$. Then, using \eqref{non-autonomous system}, it can be shown that the order conditions for \eqref{non-autonomous y(t)} are the same as those for \eqref{autonomous y(t)}. We now prove the main theorem of this subsection.

\begin{theorem}\label{thm:sl-imex-order}
Fix $\varepsilon>0$ and consider the IMEX-RK coefficients $(A,b)$ and $(\widetilde A,\widetilde b)$ that are CK and GSA. Assume that the coefficients satisfy the IMEX-RK order conditions through order $p$ for \eqref{autonomous y(t)}. Then when the same cofficients are applied to IMEX-RK scheme \eqref{eq:IMEX-RK}, the method also has local truncation error of order $p+1$, i.e. 
\begin{equation}\label{eq:sl-imex-local-order}
\|f^{n+1}-f(t_n+h)\|_{\mathcal{H}}\le C_\varepsilon h^{p+1}.
\end{equation}
where $\mathcal{H}=L^2(\mathbb{T}^d_{x}\times \mathbb{R}^d_{v})$. 
\end{theorem}

\begin{proof}
The key idea of the proof is performing a change of variable to obtain an equation of the form \eqref{non-autonomous y(t)}. Fix $\varepsilon>0$. To simplify the notations, we let  $\mathcal F_\varepsilon(f)=\frac{1}{\varepsilon}G_P(f),~\mathcal H_\varepsilon(f)=\frac{1}{\varepsilon}Q_P(f),~\mathcal{T}=-v\cdot \nabla_x$ and $(S_{\tau}z)(x,v)=z(x-\tau v,v)$. Define 
\begin{equation}\label{eq:transformed-fields-order}
\begin{aligned}
\widehat{\mathcal F}_\varepsilon(\tau,z)&=S_{-\tau}\mathcal F_\varepsilon(S_\tau z),\\
\widehat{\mathcal H}_\varepsilon(\tau,z)&=S_{-\tau}\mathcal H_\varepsilon(S_\tau z)
\end{aligned}
\end{equation}
We also define 
\begin{equation}\label{eq:integrating-factor-order}
w(\tau)=S_{-\tau}f(t_n+\tau),\qquad 0\le\tau\le h.
\end{equation}
Since
$$
\partial_t f=\mathcal T f+\mathcal F_\varepsilon(f)+\mathcal H_\varepsilon(f),
$$
differentiation of \eqref{eq:integrating-factor-order} gives
$$
w'(\tau)=-\mathcal T S_{-\tau}f(t_n+\tau)+S_{-\tau}\partial_t f(t_n+\tau) \notag=S_{-\tau}\mathcal F_\varepsilon(S_\tau w(\tau))+S_{-\tau}\mathcal H_\varepsilon(S_\tau w(\tau)).
$$
Consequently,
\begin{equation}\label{eq:transformed-ode-order}
w'(\tau)=\widehat{\mathcal F}_\varepsilon(\tau,w(\tau))+\widehat{\mathcal H}_\varepsilon(\tau,w(\tau)),\qquad w(0)=f(t_n).
\end{equation}
We therefore have obtained a non-autonomous system of the form \eqref{non-autonomous y(t)}. From the arguments there, when the IMEX-RK method with coefficients $(A,b)$ and $(\widetilde{A},\widetilde{b})$ are applied to \eqref{eq:transformed-ode-order}, we have 
\begin{equation}\label{eq:transformed-local-error-order}
\|W_{s}-w(h)\|_{\mathcal{H}} \le C_\varepsilon h^{p+1}.
\end{equation}
Moreover, the stages of the scheme are
\begin{equation}\label{eq:transformed-stages-order}
W_i=f^n+h\sum_{j=1}^{i-1}a_{ij}\widehat{\mathcal F}_\varepsilon(c_jh,W_j)+h\sum_{j=1}^{i}\widetilde a_{ij}\widehat{\mathcal H}_\varepsilon(c_jh,W_j),
\end{equation}

Define $f^{(i)}=S_{c_i h}W_i.$ Using the group property $S_\alpha S_\gamma=S_{\alpha+\gamma},$ we obtain
$$
S_{c_i h}\widehat{\mathcal F}_\varepsilon(c_jh,W_j)=S_{c_i h}S_{-c_jh}
\mathcal F_\varepsilon(S_{c_jh}W_j)=S_{(c_i-c_j)h}\mathcal F_\varepsilon(f^{(j)}).
$$
The same identity holds for $\widehat{\mathcal H}_\varepsilon$.
Multiplying \eqref{eq:transformed-stages-order} by $S_{c_i h}$
therefore gives exactly \eqref{eq:IMEX-RK}. 

Since $f(t_n+h)=S_hw(h),~ f^{n+1}=S_hW_{s}$, using the fact that $S_h$ preserves $\| \cdot \|_{\mathcal{H}}$, we have
$$
\|f^{n+1}-f(t_n+h)\|_{\mathcal{H}}=\|S_h(W_{s}-w(h))\|_{\mathcal{H}}=\|W_{s}-w(h)\|_{\mathcal{H}} \notag\le C_\varepsilon h^{p+1}.
$$
This proves the local error estimate \eqref{eq:sl-imex-local-order}.
\end{proof}

\begin{remark}
The constant $C_{\varepsilon}$ in the error estimates of \eqref{eq:sl-imex-local-order} may depend on $\varepsilon$ and need not remain bounded as
$\varepsilon\to0$. The theorem therefore establishes classical
accuracy for fixed $\varepsilon$, not uniform accuracy or the
accuracy of the limiting scheme. In fact, if $C_{\varepsilon}$ does not depend on $\varepsilon$, then the same order conditions will hold for the limiting scheme, which is typically not the case as discussed in the following subsection. 
\end{remark}

\subsubsection{Order conditions for the limiting scheme}
In this subsection, we analyze the accuracy of the limiting scheme \eqref{eq:limiting-scheme}. As mentioned previously, a reduction in accuracy order will typically be observed when the standard IMEX-RK coefficients are applied. However, we shall prove below that the macroscopic component of \eqref{eq:limiting-scheme} is at least formally first-order temporal consistent with the compressible Euler equations. For higher order accuracy, more conditions on the IMEX-RK coefficients should be imposed. We provide sufficient conditions on the coefficients to ensure second-order consistency and use it to demonstrate the additional efforts needed when higher orders are considered. Throughout the subsection, we let $h=\Delta t,~\phi(v)=\left(1,v,\frac{|v|^2}{2}\right)^T,~Pz=\int_{\mathbb R^d}\phi(v)z(v)\,dv$, and $\mathcal H_U=\mathcal L_U+\beta I$ where $\beta>0$ is constant and $\mathcal L_U$ is given by \eqref{eq:local-linearization}.

The compressible Euler equations are given by \eqref{eq:EulerEquation}, and we restate it as follows 
\begin{equation}\label{eq:ck-gsa-Euler-equations}
\partial_t U+\mathcal{T}(U)=0, 
\end{equation}
where $\mathcal{T}(U)=\nabla_x \cdot \int_{\mathbb R^d}\phi(v)\otimes v\,M[U](v)\,dv$. We first prove the following theorem, which states that the limiting scheme \eqref{eq:limiting-scheme} is formally at least first-order temporal consistent with the compressible Euler equations.   

\begin{theorem}\label{prop:ck-gsa-first-order-consistency}
Consider the CK-GSA limiting scheme \eqref{eq:limiting-scheme}
with $c_s=1$ on the periodic spatial domain $\mathbb T^d$. Suppose that
$Pg^n=0$ and write $k^{(j)}=\mathcal H_{U^{(j)}}(g^{(j)})$. If $U^n=U(t_n)$ for a smooth Euler solution $U$, then formally we have 
\begin{equation}\label{eq:ck-gsa-local-defect}
\|U^{n+1}-U(t_n+h)\|_{L_x^2}\le Ch^2.
\end{equation}
Thus, the macroscopic component of the limiting scheme is at least
first-order consistent in time with the Euler equations.
\end{theorem}

\begin{proof}
We first verify the zero-moment properties of the microscopic stage functions $g^{(j)}$. By assumption, $Pg^{(1)}=Pg^n=0$. For $j\ge2$, \eqref{eq:limiting-scheme} gives
$$
Pg^{(j)}=\frac{PR^{(j)}-PM[U^{(j)}]}{h\beta\widetilde a_{jj}}=
\frac{U^{(j)}-U^{(j)}}{h\beta\widetilde a_{jj}}=0.
$$
Moreover, collision conservation implies $P\mathcal L_{U^{(j)}}(g^{(j)})=0,$
and therefore
\begin{equation}\label{eq:ck-gsa-zero-moments}
Pk^{(j)}=P\mathcal H_{U^{(j)}}g^{(j)}=P\mathcal L_{U^{(j)}}(g^{(j)})+\beta Pg^{(j)}=0.
\end{equation}
Taking Taylor expansion with respect to the space variable $x$ gives 
\begin{equation}\label{Taylor expansion for g,k wrt x}
\begin{aligned}
& S_{i,j}[g^{(j)}]=g^{(j)}(x-(c_i - c_j)hv,v)=g^{(j)}(x,v)-(c_i-c_j)hv\cdot\nabla_x g^{(j)}(x,v)+O(h^2). \\
& S_{i,j}[k^{(j)}]=k^{(j)}(x-(c_i - c_j)hv,v)=k^{(j)}(x,v)-(c_i-c_j)hv\cdot\nabla_x k^{(j)}(x,v)+O(h^2)
\end{aligned}
\end{equation}

Taking moments on both sides and use the fact that $Pg^{(j)}=Pk^{(j)}=0$ gives  \begin{equation}\label{eq:ck-gsa-stage-moment-bounds}
PS_{i,j}[g^{(j)}]=O(h),\qquad PS_{i,j}[k^{(j)}]=O(h), \qquad j<i.
\end{equation}

Taking moments in \eqref{eq:limiting-R}, we therefore formally obtain
\begin{equation}\label{eq:ck-gsa-leading-macro-stage}
\begin{aligned}
U^{(i)}={}&PS_{i,0}[M[U^n]]+h\sum_{j=1}^{i-1}a_{ij}PS_{i,j}[k^{(j)}]-h\beta\sum_{j=1}^{i-1}\widetilde a_{ij}PS_{i,j}[g^{(j)}]\\
={}&PS_{i,0}[M[U^n]]+O(h^2).
\end{aligned}
\end{equation}
Taylor expansion of the shifted Maxwellian gives
$$
S_{i,0}[M[U^n]]=M[U^n]-c_i h\,v\cdot\nabla_x M[U^n]+O(h^2).
$$
Using $PM[U^n]=U^n$ and $P\bigl(v\cdot\nabla_x M[U^n]\bigr)= \mathcal{T}(U^n)$, we obtain
\begin{equation}\label{Taylor expand of PSM}
PS_{i,0}[M[U^n]]=U^n-c_i h\mathcal{T}(U^n)+O(h^2).
\end{equation}

Finally, GSA condition and $c_s=1$ give
\begin{equation}\label{final Taylor expansion 1}
U^{n+1}=U^{(s)}=U^n-h\mathcal{T}(U^n)+O(h^2).
\end{equation}

For a smooth solution of \eqref{eq:ck-gsa-Euler-equations},
Taylor expansion in time yields
\begin{equation}\label{final Taylor expansion 2}
U(t_n+h)=U(t_n)+h\partial_tU(t_n)+O(h^2)=U(t_n)-h\mathcal{T}(U)(t_n)+O(h^2).
\end{equation}
Taking $U^n=U(t_n)$ and comparing the two expansions \eqref{final Taylor expansion 1}, \eqref{final Taylor expansion 2} proves \eqref{eq:ck-gsa-local-defect}.
\end{proof}

We now analyze the second-order accuracy of the limiting scheme. We continue to consider CK and GSA schemes with $c_s=1$. We assume the standard second-order conditions for \eqref{eq:IMEX-RK}(or \eqref{autonomous y(t)}) given by 
\begin{equation}\label{eq:second-classical-conditions}
b^T\mathbf 1=\widetilde b^T\mathbf 1=1,\qquad b^Tc=\widetilde b^Tc=\frac12.
\end{equation}

To analyze the second-order temporal accuracy of \eqref{eq:limiting-scheme}, we need to study the terms in the Taylor expansion of $U^{n+1}$ up to $h^2$. To this end, we first notice that from \eqref{Taylor expand of PSM}, we have $U^{(i)}=U^n-c_i h\mathcal{T}(U^n)+O(h^2)$. Hence, $M[U^{(i)}]=M[U^n]+O(h)$. From the definition of $R^{(i)}$ in \eqref{eq:limiting-R}, it is easy to see that $R^{(i)}=M[U^n]+O(h)$. Therefore, $R^{(i)}-M[U^{(i)}]=O(h)$ and from the definition of $g^{(i)}$ in \eqref{eq:limiting-scheme}, we have $g^{(i)}=O(1)$ with respect to $h$ and hence 
\begin{equation}\label{eq:second-g-stage-leading}
g^{(i)}=\gamma_i+O(h),\qquad\gamma_1=g^n.
\end{equation}
To simplify the notations, we define the following
\begin{equation}\label{def of D,T,Z}
D=v\cdot\nabla_x,\qquad  \mathcal{T}(U)=PDM[U],\qquad Z=\frac{\mathcal{H}_U}{\beta}=I+\frac{\mathcal L_U}{\beta}.
\end{equation}
Then, we can prove the following lemma. 

\begin{lemma}\label{second order expansion for U^n+1}
Consider the CK--GSA limiting scheme \eqref{eq:limiting-scheme} with $c_s=1$, then 
\begin{equation}\label{eq:second-macro-expansion-last stage}
\begin{aligned}
U^{n+1}=U^{(s)}&=U^n-h\mathcal{T}(U^n)+\frac{h^2}{2}PD^2M[U^n]\\
&\quad{}+h^2\beta PD\sum_{j<s}(1-c_j)(\widetilde a_{sj}I-a_{sj}Z)\gamma_j+O(h^3).
\end{aligned}
\end{equation}
\end{lemma}

\begin{proof}
From \eqref{eq:ck-gsa-leading-macro-stage}, we have 
$$
U^{(i)}=PS_{i,0}[M[U^n]]+h\sum_{j<i}a_{ij}P S_{i,j}[k^{(j)}]-h\beta\sum_{j<i}\widetilde a_{ij}P S_{i,j}[g^{(j)}].
$$

Same as the 1st order case, we still have $Pg^{(i)}=Pk^{(i)}=0$. To analyze the second-order accuracy, we need to perform Taylor expansion on $PS_{i,0}[M[U^n]]$ up to $h^2$ and on $P S_{i,j}[k^{(j)}], P S_{i,j}[g^{(j)}]$ up to $h$. Similar to \eqref{Taylor expand of PSM}, performing Taylor expansion on $PS_{s,0}[M[U^n]]$ with respect to $x$ up to $h^2$ gives 

\begin{equation}\label{eq:PS(MU)}
\begin{aligned}
PS_{i,0}[M[U^n]]=&P(M[U^n]-c_ihDM[U^n]+\frac{c_i^{2}h^2}{2}D^2M[U^n]+O(h^3)) \\
=&U^n - c_ih\mathcal{T}(U^n)+\frac{c_i^{2}h^2}{2}PD^2M[U^n]+O(h^3)
\end{aligned}
\end{equation}

We now analyze the situation for $P S_{i,j}[k^{(j)}], P S_{i,j}[g^{(j)}]$. From \eqref{Taylor expansion for g,k wrt x} and \eqref{eq:second-g-stage-leading}, we have  
$$
\begin{aligned}
P S_{i,j}[g^{(j)}]&=-(c_i-c_j)h\,PD\gamma_j+O(h^2),\\
P S_{i,j}[k^{(j)}]&=-(c_i-c_j)h\,PD(\beta Z\gamma_j)+O(h^2),
\end{aligned}
$$
Consequently,
\begin{equation}\label{eq:second-macro-expansion}
\begin{aligned}
U^{(i)}&=U^n-c_ih\mathcal{T}(U^n)+\frac{c_i^2h^2}{2}PD^2M[U^n]\\
&\quad{}+h^2\beta PD\sum_{j<i}(c_i-c_j)(\widetilde a_{ij}I-a_{ij}Z)\gamma_j+O(h^3).
\end{aligned}
\end{equation}
In particular, for $c_s=1$, we have $U^{n+1}=U^{(s)}$ and \eqref{eq:second-macro-expansion-last stage} is obtained. 
\end{proof}

The smooth Euler solution starting from $t_n$ satisfies
\begin{equation}\label{eq:second-exact-Euler-expansion}
U(t_n+h)=U(t_n)-h\mathcal{T}(U)(t_n)+\frac{h^2}{2}\mathcal{T}'(U)\mathcal{T}(U)(t_n)+O(h^3).
\end{equation}
To derive the conditions for second-order accuracy, we need to consider the term containing $h^2$ in $U^{n+1}-U(t_n+h)$ and simplify the final expression. To this end, we define polynomials $p_i(z)$ and $q_i(z)$ recursively by $p_1(z)=1,~q_1(z)=0$, and for $i=2,\ldots,s$,
\begin{equation}\label{eq:second-p-q-recursion}
\widetilde a_{ii}p_i(z)=\sum_{j<i}(a_{ij}z-\widetilde a_{ij})p_j(z),\qquad  \widetilde a_{ii}q_i(z)=c_i+\sum_{j<i}(a_{ij}z-\widetilde a_{ij})q_j(z).
\end{equation}
We also define
\begin{equation}\label{eq:second-rho-sigma}
\rho(z)=\sum_{j<s}(1-c_j)(\widetilde a_{sj}-a_{sj}z)p_j(z),\qquad \sigma(z)=\sum_{j<s}(1-c_j)(\widetilde a_{sj}-a_{sj}z)q_j(z).
\end{equation}

We view $M[\cdot]$ as an operator from the space of macroscopic states to some function space, and let $J_{U^n}$ be the Fr\'echet derivative of $M[\cdot]$ at $U^n$. In particular, $J_{U^n}$ is a linear operator with the same domain and target space as $M[\cdot]$. Then we let 
\begin{equation}\label{def of Pi and W}
\Pi_{U^n}=J_{U^n}P, \qquad W=(I-\Pi_{U^n})DM[U^n]
\end{equation}
Using the polynomials in \eqref{eq:second-rho-sigma} and the variables in \eqref{def of Pi and W}, we can prove the following:

\begin{lemma}\label{second order defect in h}
Consider the CK--GSA limiting scheme \eqref{eq:limiting-scheme} with $c_s=1$, then 
\begin{equation}\label{eq:second-defect-polynomial}
U^{n+1}-U(t_n+h)=h^2PD\left[\left(\frac12I- \sigma(Z)\right)W+\beta\rho(Z)g^n\right]+O(h^3).
\end{equation}
\end{lemma}

\begin{proof}
Since $PM[U]=U$ for any $U$, taking differentiation on both sides at $U^n$ and using the linearity of $P$ gives $PJ_{U^n}=I$. Consequently,
\begin{equation}\label{eq:DM=J+W}
DM[U^n]=\Pi_{U^n}DM[U^n]+(I-\Pi_{U^n})DM[U^n]=J_{U^n}\mathcal{T}(U^n)+W,\qquad PW=0.
\end{equation}
Notice that $\mathcal{T}(U)=PDM[U]$ for any $U$. Taking differention on both sides at $U^n$ and using the linearity of $P$ and $D$ gives $\mathcal{T}'(U^n)=PDJ_{U^n}$. Therefore, from \eqref{eq:DM=J+W}, we have 
\begin{equation}\label{eq:second-Euler-identity}
PD^2M[U^n]-\mathcal{T}'(U^n)\mathcal{T}(U^n)=PDW.
\end{equation}
Using \eqref{eq:second-macro-expansion-last stage},\eqref{eq:second-exact-Euler-expansion} and \eqref{eq:second-Euler-identity}, the coefficient of the
second-order one-step defect in $U^{n+1}-U(t_n+h)$ is therefore
\begin{equation}\label{eq:second-defect-gamma}
U^{n+1}-U(t_n+h)=h^2PD\left[\frac12W+\beta\sum_{j<s}(1-c_j)(\widetilde a_{sj}I-a_{sj}Z)\gamma_j\right]+O(h^3).
\end{equation}

Now, we want to describe $\gamma_j$ in terms of $M[U^n]$ and $g^n$. Recall that the limiting scheme updates \eqref{eq:limiting-scheme} can be written as
\begin{equation}\label{rewrite limiting scheme}
M[U^{(i)}]+h\beta\widetilde a_{ii}g^{(i)}=S_{i,0}[M[U^n]]+h\sum_{j<i}S_{i,j}
\left[a_{ij}k^{(j)}-\beta\widetilde a_{ij}g^{(j)}\right].
\end{equation}
We want to expand \eqref{rewrite limiting scheme} to first order in $h$ on both sides. Notice that taking Taylor expansion of $M[\cdot]$ at $U^n$ from \eqref{eq:second-macro-expansion} gives 
\begin{equation}\label{Taylor expansion of M at Un}
M[U^{(i)}]=M[U^n]-c_i hJ_U\mathcal{T}(U^n)+O(h^2).
\end{equation}
Then from \eqref{Taylor expansion for g,k wrt x}, \eqref{eq:PS(MU)}, \eqref{eq:second-g-stage-leading},\eqref{Taylor expansion of M at Un} and $\mathcal H_U=\beta Z$, we can expand \eqref{rewrite limiting scheme} and compare the terms involving first-order $h$. This yields
$$
-c_iJ_{U^n}\mathcal{T}(U^n)+\beta\widetilde a_{ii}\gamma_i=-c_iDM[U^n]+\beta\sum_{j<i}(a_{ij}Z-\widetilde a_{ij}I)\gamma_j.
$$
Therefore,
\begin{equation}\label{eq:second-gamma-recursion}
\widetilde a_{ii}\gamma_i=-\frac{c_i}{\beta}W+\sum_{j<i}(a_{ij}Z-\widetilde a_{ij}I)\gamma_j.
\end{equation}
Equation \eqref{eq:second-gamma-recursion}, together with a simple induction argument, implies
\begin{equation}\label{eq:second-gamma-polynomial}
\gamma_i=p_i(Z)g^n-\frac1\beta q_i(Z)W.
\end{equation}
Substitution of \eqref{eq:second-gamma-polynomial} into \eqref{eq:second-defect-gamma} gives \eqref{eq:second-defect-polynomial} as desired. 
\end{proof}

Now, we are ready to state and prove the suffcient conditions for the second-order temporal consistency of the limiting scheme \eqref{eq:limiting-scheme}. 

\begin{theorem}\label{thm:ck-gsa-second-order-consistency}
Consider the CK--GSA limiting scheme \eqref{eq:limiting-scheme}, with $c_s=1$ and the standard second-order conditions \eqref{eq:second-classical-conditions}.
Let $p_i(z)$ be defined by \eqref{eq:second-p-q-recursion}. If the additional coefficient identity
\begin{equation}\label{eq:second-additional-condition}
\rho(z) = \sum_{j<s}(1-c_j)(\widetilde a_{sj}-a_{sj}z)p_j(z)\equiv0
\end{equation}
holds, then
$$
\|U^{n+1}-U(t_n+h)\|_{L_{x}^2}\le Ch^3.
$$
\end{theorem}

\begin{proof}
From \eqref{eq:second-defect-polynomial}, to ensure second-order temporal consistency, we need the $h^2$ term to vanish identically, i.e. $\rho(z)\equiv0,~\sigma(z)\equiv\frac12$. Therefore, it remains to prove that when $\rho(z) \equiv 0$, the condition $\sigma(z)\equiv\frac12$ is automatically satisfied. 

It can be easily checked by a simple induction argument that the following useful identity for $p_i,q_i$ holds: 
\begin{equation}\label{eq:second-pq-identity}
p_i(z)+(1-z)q_i(z)=1.
\end{equation}

Using \eqref{eq:second-pq-identity} and the standard second-order RK conditions \eqref{eq:second-classical-conditions}, we have 
$$
\rho(z)+(1-z)\sigma(z)=\sum_{j<s}(1-c_j)(\widetilde a_{sj}-a_{sj}z)=\widetilde b^T(\mathbf1-c)-zb^T(\mathbf1-c)=\frac{1-z}{2}.
$$
Hence, if $\rho(z)\equiv 0$, this polynomial identity implies $\sigma(z)\equiv \frac{1}{2}$.
\end{proof}

\subsection{Stability analysis}\label{sec:stability}

In this section, we analyze the stability of the scheme \eqref{eq:IMEX-RK} and the limiting scheme \eqref{eq:limiting-scheme}. The key ingredient of the stability analysis is to construct an $\varepsilon$-dependent total energy(see \eqref{eq:euler-energy}) which converges to a proper total energy in the limiting scheme(\eqref{eq:euler-limit-energy}) as $\varepsilon \to 0$. We prove that both the $\varepsilon$-dependent energy and the limiting energy remain nonincreasing over time, and therefore establish the stability of our scheme across different regimes.   

We emphasize that the stability analysis presented here concerns the linearization of the proposed numerical scheme \eqref{eq:IMEX-RK} about a fixed global Maxwellian. More precisely, we establish an energy estimate for the discrete first-order perturbation dynamics, rather than a nonlinear stability theorem for the full scheme. This analysis identifies how the collision structure, characteristic transport, and IMEX stage coupling jointly control perturbation growth, and provides a foundation for a subsequent nonlinear analysis. Extending the result to the full nonlinear scheme requires additional estimates for the nonlinear collision operator, which we leave to a future study. We also comment that in this section, we only work on the two-stage IMEX Euler method to demonstrate our main methodology. 

We first introduce the setting of the linearziation problem and present the linearized scheme. Then we establish a complete stability analysis for the two-stage IMEX Euler method of \eqref{eq:IMEX-RK}. Finally, we estabilish similar stability property of the limiting scheme \eqref{eq:limiting-scheme} with the same IMEX Euler coefficients. Throughout this subsection, $\varepsilon>0$ is fixed during each computation, and the time step $h=\Delta t>0$ is constant. We use the subscript $\varepsilon$ to emphasize that the variables in consideration depend on the fixed choice of $\varepsilon$.

\subsubsection{Setting and linearization}
Let
\begin{equation}
M_*(v)=\frac{\rho_*}{(2\pi\theta_*)^{d/2}}\exp\left(-\frac{|v-u_*|^2}{2\theta_*}\right),\qquad \rho_*,\theta_*>0,
\end{equation}
be a fixed global Maxwellian. A natural choice of $M_*$ is the Maxwellian with the same total mass, momentum, and energy as the initial distribution. We introduce the normalized perturbation around $M_{*}$
\begin{equation}\label{eq:stability-perturbation}
f_{\varepsilon}=M_*+\delta\sqrt{M_*}\,w_{\varepsilon}+O(\delta^2),
\end{equation}
where $\delta$ is independent of $\varepsilon$. Substituting \eqref{eq:stability-perturbation} into the original Boltzmann equation \eqref{eq:BoltzmannEquation} and remove the higher order $\delta$ terms yields the linearized Boltzmann equation: 
\begin{equation}\label{linearized Boltz eq}
\partial_t w_{\varepsilon} + v \cdot \nabla_{x} w_{\varepsilon} = \frac{1}{\mathrm{\varepsilon}} L(w_{\varepsilon}),
\end{equation}
where $L(w_{\varepsilon})=M_*^{-1/2}\left[Q(M_*,\sqrt{M_*}w_{\varepsilon})+Q(\sqrt{M_*}w_{\varepsilon},M_*)\right]$. Let $\Pi$ denote the orthogonal projection in $L^2_v$ onto $\mathcal N:=\ker L$, where
\begin{equation}\label{kernel space of L}
\mathcal N=\operatorname{span}\left\{\sqrt{M_*},(v_j-u_{*,j})\sqrt{M_*},(|v-u_*|^2-d\theta_*)\sqrt{M_*}\right\}.
\end{equation}
Write $\Pi^\perp=I-\Pi$. We introduce the macroscopic and scaled microscopic variables
\begin{equation}\label{eq:stability-pq}
p_{\varepsilon}=\Pi w_{\varepsilon},\qquad q_{\varepsilon}=\varepsilon^{-1}\Pi^\perp w_{\varepsilon},\qquad w_{\varepsilon}=p_{\varepsilon}+\varepsilon q_{\varepsilon}.
\end{equation}

\begin{remark}
If the original local decomposition is $f_{\varepsilon}=M[U_{\varepsilon}]+\varepsilon g_{\varepsilon}$, then at first order in $\delta$,
\begin{equation}
M[U_{\varepsilon}]=M_*+\delta\sqrt{M_*}p_{\varepsilon}+O(\delta^2),\qquad g_{\varepsilon}=\delta\sqrt{M_*}q_{\varepsilon}+O(\delta^2),
\end{equation}
where these expansions are understood at fixed $\varepsilon$.
\end{remark}

Now, we assume that $L$ is bounded, self-adjoint, and nonpositive, with the spectrum range given by  
\begin{equation}\label{eq:stability-spectrum}
\sigma(L|_{\mathcal N^\perp})\subset[-\Lambda,-\lambda],\qquad 0<\lambda\le\Lambda<\infty.
\end{equation}
Equivalently, on the microscopic subspace $\mathcal{N}^\perp$, we have 
\begin{equation}
-\Lambda\|z\|_{L_v^2}^2\le\langle Lz,z\rangle_{L_v^2}\le-\lambda\|z\|_{L_v^2}^2.
\end{equation}
Note that assumption \eqref{eq:stability-spectrum} holds, in particular, for Maxwell molecules.

\subsubsection{Stability for the two-stage IMEX Euler method}
We now analyze the stability of the IMEX Euler scheme. Consider scheme \eqref{eq:IMEX-RK} with 
\begin{equation}\label{cofficients for IMEX euler}
A=\begin{pmatrix}
    0&0\\
    1&0
\end{pmatrix},\qquad
\widetilde A=\begin{pmatrix}
    0&0\\
    0&1
\end{pmatrix},\qquad
c=\begin{pmatrix}0\\1\end{pmatrix}.
\end{equation}
We further assume the scheme is GSA so that $b^T=(1,0), \widetilde b^T=(0,1)$. Then the scheme update can be written explicitly as  
\begin{equation}\label{update for IMEX Euler}
f^{n+1}_{\varepsilon}=\mathcal{S}_hf^n_{\varepsilon}+\frac{h}{\varepsilon}\mathcal{S}_hG_P(f^n_{\varepsilon})+\frac{h}{\varepsilon}Q_P(f^{n+1}_{\varepsilon}),
\end{equation}
where we write $\mathcal{S}_{h}=\mathcal{S}_{v,\Delta t}$ to simplify the notation. \eqref{update for IMEX Euler} is a scheme for the original Boltzmann equation. Since we are now interested in the linearized equation \eqref{linearized Boltz eq}, we shall first derive the linearization of the scheme \eqref{update for IMEX Euler} to obtain an update for the first variation $w_{\varepsilon}$. 

The key of the derivation is to compute the Fr\'echet derivatives of the operators $G_P$ and $Q_P$ at the global Maxwellian $M_*$, which we denote by $DQ_P(M_*),~ DQ_P(M_*)$ repspectively. It can be computed directly that the derivative of the Maxwellian reconstruction is
\begin{equation}
D\mathcal{M}(M_*)[\sqrt{M_*}w_{\varepsilon}]=\sqrt{M_*}\Pi w_{\varepsilon}.
\end{equation}
Consequently,
\begin{equation}\label{eq:stability-linearized-collision}
\begin{aligned}
M_*^{-1/2}DQ_P(M_*)[\sqrt{M_*}w_{\varepsilon}]&=-\beta\Pi^\perp w_{\varepsilon},\\
M_*^{-1/2}DG_P(M_*)[\sqrt{M_*}w_{\varepsilon}]&=(L+\beta\Pi^\perp)w_{\varepsilon}.
\end{aligned}
\end{equation}
Using \eqref{eq:stability-linearized-collision}, we obtain
\begin{equation}\label{linear scheme for IMEX euler}
w^{n+1}_{\varepsilon}=S_hw^n_{\varepsilon}+\frac{h}{\varepsilon}S_h(L+\beta\Pi^\perp)w^n_{\varepsilon}-\frac{h\beta}{\varepsilon}\Pi^\perp w^{n+1}_{\varepsilon}.
\end{equation}
This is precisely the linearized scheme of \eqref{update for IMEX Euler} we are interested in. Set $r=\frac{h\beta}{\varepsilon}$ and define 
\begin{equation}\label{def of C}
C=\left.\left(I+\frac{L}{\beta}\right)\right|_{\mathcal N^\perp}.
\end{equation}
Then $L+\beta\Pi^\perp=\beta C\Pi^\perp$ and the update \eqref{linear scheme for IMEX euler} becomes
\begin{equation}\label{eq:euler-linear-update}
(I+r\Pi^\perp)w^{n+1}_{\varepsilon}=S_h\bigl[\Pi+(I+rC)\Pi^\perp\bigr]w^n_{\varepsilon}.
\end{equation}

Now, consider the working space $\mathcal{H}=L^2(\mathbb T_x^d\times\mathbb R_v^d)$. Define the total energy
\begin{equation}\label{eq:euler-energy}
E_{\varepsilon,h}^n:=\|p^n_{\varepsilon}\|_{\mathcal H}^2+(\varepsilon+h\beta)^2\|q^n_{\varepsilon}\|_{\mathcal H}^2.
\end{equation}
We now prove the following stability result for the linearized scheme \eqref{linear scheme for IMEX euler}.

\begin{theorem}\label{stability for linear IMEX euler}
Consider the IMEX Euler scheme \eqref{linear scheme for IMEX euler} for the linearized Boltzmann equation \eqref{linearized Boltz eq}. Assume \eqref{eq:stability-spectrum} holds for $L$ for some $0<\lambda \le \Lambda < \infty$ and assume $\beta\geq \frac{\Lambda}{2}$. Then we have  
\begin{equation}\label{eq:euler-energy-stability}
E_{\varepsilon,h}^{n+1}\le E_{\varepsilon,h}^n\le E_{\varepsilon,h}^0,\qquad n\ge0.
\end{equation}
\end{theorem}

\begin{proof}
First, notice that periodicity of the spatial domain implies
\begin{equation}\label{eq:shift isometry}
\begin{aligned}
\|S_\tau z\|_{\mathcal H}^2=\int_{\mathbb R^d}\int_{\mathbb T^d}|z(x-\tau v,v)|^2\,dx\,dv=\|z\|_{\mathcal H}^2.
\end{aligned}
\end{equation}
Moreover, $S_\tau S_\sigma=S_{\tau+\sigma}$.

Now, define $y^n_{\varepsilon}:=(I+r\Pi^\perp)w^n_{\varepsilon}=p^n_{\varepsilon}+(\varepsilon+h\beta)q^n_{\varepsilon}$. Then $E_{\varepsilon,h}^n=\|y^n_{\varepsilon}\|_{\mathcal H}^2$. Since $w^n_{\varepsilon}=\left(\Pi+\frac{1}{1+r}\Pi^\perp\right)y^n_{\varepsilon}$, equation \eqref{eq:euler-linear-update} gives $y^{n+1}_{\varepsilon}=S_h(\Pi+C_{\varepsilon,h}\Pi^\perp)y^n_{\varepsilon}$,
where $C_{\varepsilon,h}=\frac{I+rC}{1+r}=I+\frac{h}{\varepsilon+h\beta}L$ acts on $\mathcal N^\perp$. 

The operator norm below is the induced norm on the microscopic space:
\begin{equation}
\|C_{\varepsilon,h}\|_{\mathrm{mic}}:=\sup_{\substack{z\in\mathcal N^\perp\\z\ne0}}\frac{\|C_{\varepsilon,h}z\|_{L_v^2}}{\|z\|_{L_v^2}}.
\end{equation}
By self-adjointness and the spectral theorem,
\begin{equation}\label{eq:euler-spectral-norm}
\|C_{\varepsilon,h}\|_{\mathrm{mic}}=\sup_{\ell\in\sigma(L|_{\mathcal N^\perp})}\left|1+\frac{h\ell}{\varepsilon+h\beta}\right|.
\end{equation}
Since $\beta\ge\frac{\Lambda}{2}$, we have $0<\frac{h\Lambda}{\varepsilon+h\beta}\le2$. Therefore, $\|C_{\varepsilon,h}\|_{\mathrm{mic}}\le 1$. Using \eqref{eq:shift isometry}, we can derive 
\begin{equation}
E_{\varepsilon,h}^{n+1}=\|S_h(\Pi+C_{\varepsilon,h}\Pi^\perp)y^n_{\varepsilon}\|_{\mathcal H}^2=\|(\Pi+C_{\varepsilon,h}\Pi^\perp)y^n_{\varepsilon}\|_{\mathcal H}^2.
\end{equation}
The macroscopic and microscopic subspaces are orthogonal, and $C_{\varepsilon,h}$ preserves the microscopic subspace. Hence
\begin{equation}\label{eq:euler-exact-energy-identity}
E_{\varepsilon,h}^{n+1}=\|\Pi y^n_{\varepsilon}\|_{\mathcal H}^2+ \|C_{\varepsilon,h}\Pi^\perp y^n_{\varepsilon}\|_{\mathcal H}^2=\|p^n_{\varepsilon}\|_{\mathcal H}^2
+(\varepsilon+h\beta)^2\|C_{\varepsilon,h}q^n_{\varepsilon}\|_{\mathcal H}^2.
\end{equation}
Applying the fact that $\|C_{\varepsilon,h}\|_{\mathrm{mic}}\le 1$ gives 
\eqref{eq:euler-energy-stability}. This completes the proof. 
\end{proof}

\subsection{Stability for the linearized limiting scheme}
In this subsection, we prove that the stability result we derived in the previous subsection is asymptotically preserved as $\varepsilon\to 0$. We start by finding the limiting scheme of the linearized IMEX Euler scheme \eqref{linear scheme for IMEX euler}. From \eqref{eq:euler-linear-update} and using $w^n_{\varepsilon}=p^n_{\varepsilon}+\varepsilon q^n_{\varepsilon}$,
we can obtain
\begin{equation}\label{eq:euler-linear-micro-macro}
p^{n+1}_{\varepsilon}+(\varepsilon+h\beta)q^{n+1}_{\varepsilon}=S_h\left[p^n_{\varepsilon}+\varepsilon q^n_{\varepsilon}+h\beta Cq^n_{\varepsilon}\right].
\end{equation}
Applying $\Pi$ and $\Pi^\perp$, respectively, gives
\begin{equation}\label{eq:euler-linear-micro-macro-projected}
\begin{aligned}
p^{n+1}_{\varepsilon}&=\Pi S_h\left[p^n_{\varepsilon}+\varepsilon q^n_{\varepsilon}+h\beta Cq^n_{\varepsilon}\right],\\
q^{n+1}_{\varepsilon}&=\frac{1}{\varepsilon+h\beta}\Pi^\perp S_h\left[p^n_{\varepsilon}+\varepsilon q^n_{\varepsilon}+h\beta Cq^n_{\varepsilon}\right].
\end{aligned}
\end{equation}
Suppose that the initial data satisfy $p_\varepsilon^0\to p^0,~ q_\varepsilon^0\to q^0$ in $\mathcal{H}$ as $\varepsilon \to0$. Under proper boundedness assumptions on the operators, it can be shown that for every fixed number of time steps $n$, $p_\varepsilon^n\to p^n$ and $q_\varepsilon^n\to q^n$. Thus, for fixed $h>0$, letting $\varepsilon\to 0$ yields the coupled limiting scheme
\begin{equation}\label{eq:euler-linear-limit}
\begin{aligned}
R^n&=S_h\left[p^n+h\beta Cq^n\right],\\
p^{n+1}&=\Pi R^n,\\
q^{n+1}&=\frac{1}{h\beta}\Pi^\perp R^n,
\end{aligned}
\end{equation}
where $C$ is defined in \eqref{def of C} and $p^n\in\mathcal N,~ q^n\in\mathcal N^\perp$. Equivalently,
\begin{equation}\label{eq:euler-linear-limit-combined}
p^{n+1}+h\beta q^{n+1}=S_h\left[p^n+h\beta Cq^n\right].
\end{equation}

This scheme is precisely the normalized linearization around $M_*$ of the nonlinear coupled limiting scheme \eqref{eq:limiting-scheme}-\eqref{eq:limiting-R}, where the coefficient matrices are chosen to be \eqref{cofficients for IMEX euler}. The Maxwellian perturbation corresponds to $p$ while the microscopic correction corresponds to $q$. Also, as $\varepsilon \to 0$, $E_{\varepsilon,h}^n\longrightarrow E_{0,h}^n$, where $E_{\varepsilon,h}^n$ is given in \eqref{eq:euler-energy} and $E_{0,h}^n$ is the energy defined by 
\begin{equation}\label{eq:euler-limit-energy}
E_{0,h}^n:=\|p^n\|_{\mathcal H}^2+h^2\beta^2\|q^n\|_{\mathcal H}^2.
\end{equation}
We can prove the following theorem, which states that the limiting scheme \eqref{eq:euler-linear-limit} is stable with respect to this limiting energy. 

\begin{theorem}\label{thm:euler-limit-stability}
Consider the limiting scheme \eqref{eq:euler-linear-limit}. Assume \eqref{eq:stability-spectrum} holds for $L$ for some $0<\lambda \le \Lambda < \infty$ and assume $\beta\geq \frac{\Lambda}{2}$. Then, for every fixed $h>0$, we have 
\begin{equation}\label{eq:euler-limit-energy-stability}
E_{0,h}^{n+1}\le E_{0,h}^n\le E_{0,h}^0,\qquad n\ge0.
\end{equation}
\end{theorem}

\begin{proof}
The proof is completely parallel to the proof of Theorem \ref{stability for linear IMEX euler}. Define $y^n:=p^n+h\beta q^n$. By orthogonality of the macroscopic and microscopic subspaces, $E_{0,h}^n=\|y^n\|_{\mathcal H}^2$,
and $\Pi y^n=p^n,~\Pi^\perp y^n=h\beta q^n.$ Consequently, \eqref{eq:euler-linear-limit-combined} becomes
\begin{equation}\label{eq:euler-limit-y-update}
y^{n+1}=S_h(\Pi+C\Pi^\perp)y^n.
\end{equation}

By self-adjointness and the spectral theorem,
\begin{equation}\label{eq:euler-limit-spectral-norm}
\|C\|_{\mathrm{mic}}=\sup_{\ell\in\sigma(L|_{\mathcal N^\perp})}\left|1+\frac{\ell}{\beta}\right|.
\end{equation}
Since $\ell\in[-\Lambda,-\lambda], ~\beta\ge\frac{\Lambda}{2}$, we have $-1\le1+\frac{\ell}{\beta}\le1$. Hence, $\|C\|_{\mathrm{mic}}\le1$.

Using the isometry of the shift operator on the periodic domain, we obtain
\begin{equation}
E_{0,h}^{n+1}=\|S_h(\Pi+C\Pi^\perp)y^n\|_{\mathcal H}^2=\|(\Pi+C\Pi^\perp)y^n\|_{\mathcal H}^2.
\end{equation}
Since $C$ preserves the microscopic subspace, orthogonality gives the
exact identity
\begin{equation}\label{eq:euler-limit-exact-energy}
E_{0,h}^{n+1}=\|p^n\|_{\mathcal H}^2+h^2\beta^2\|Cq^n\|_{\mathcal H}^2.
\end{equation}
Therefore,
\begin{equation}
E_{0,h}^{n+1}\le\|p^n\|_{\mathcal H}^2+h^2\beta^2\|q^n\|_{\mathcal H}^2=E_{0,h}^n.
\end{equation}
Iteration yields \eqref{eq:euler-limit-energy-stability}.
\end{proof}

\begin{remark}
    When the real solution is far from a global equilibrium, the linear approximation is no longer valid and a more delicate analysis is required to establish the stability of the scheme. We defer this to a future study. 
\end{remark}

\section{Numerical examples}
\label{sec:numerics}
In this section, we apply the nodal AP-IMEX-SLDG methods described in section~\ref{sec:SLDG} to several 1D2V problems to evaluate the effectiveness and efficiency of our proposed scheme. In our tests, we make the Maxwell molecule assumption that corresponds to a constant collision kernel in~\eqref{eq:BoltzmannOperator}, with $B = 1/2\pi$. The collision operator $Q$ is computed using the fast spectral method in~\cite{Mouhot-Pareschi}.
The spatial domain \( x \in [0, 1] \) is partitioned into $N_x$ elements. 
For the velocity discretization, we consider the velocity domain \( v \in [-7, 7] \)  with \( N_v = 32 \) velocity points. To test the performance of our schemes under different CFL conditions, we employ a small CFL value of $0.5$ and a larger value of 2 in all experiments. We define the time step by $\Delta t = \text{CFL} \frac{\Delta x}{v_{\text{max}}}$.

\subsection*{Test I: Accuracy test}
We first study the accuracy order of the proposed schemes. The Maxwellian distribution
{\small
\begin{equation*}
    f (t=0, x, v) = \frac{\rho(x)}{2\pi T(x)} \exp \left(- \frac{(v - U(x))^2}{2T(x)}\right)
\end{equation*}
}is assumed as initial data, with density, velocity and temperature given by
{\small
\begin{equation}
\label{eq:macro_in_smooth_test}
    \rho(x) = \frac{2 + \sin (2\pi x)}{2}, \quad U(x) = (3/4, -3/4), \quad T(x) = \frac{5 + 2\cos (2\pi x)}{20}.
\end{equation}
}Periodic boundary condition is assumed. 
We set \( N_x = 4, 8, 16, 32 \) and compute the corresponding relative $\ell^1$ and $\ell^2$ errors of the density of the solution at final time $t_{\text{final}} = 0.1$, with errors defined by
{\small
    $
    e^1_{N_x} = \|\rho_{f_{N_x}} - \rho_{f_{N_x/2}}\|_{\ell^1}/\|\rho_{f_{N_x}}\|_{\ell^1}, \; e^2_{N_x} = \|\rho_{f_{N_x}} - \rho_{f_{N_x/2}}\|_{\ell^2}/\|\rho_{f_{N_x}}\|_{\ell^2}, $
}
where $f_{N_x}$ denotes the numerical solution obtained by using $N_x$ in spatial discretization. A numerical scheme is $p$-th order accurate if $e_{N_x} \leq \frac{C}{N_x^p}$ for sufficiently large $N_x$. For the nodal-DG reconstruction, we use Lagrangian polynomials of degree $k=1$ and $k=2$.
Thus the numerical schemes are expected to achieve second and third order convergence, respectively. For the temporal discretization, the type CK method \texttt{ARS443} is used. 
\begin{table}[!h]
    \centering
    \setlength{\tabcolsep}{3pt}
    \resizebox{\textwidth}{!}{%
    \begin{tabular}{c c | cc cc | cc cc || cc cc | cc cc}
    \toprule
    & & \multicolumn{8}{c||}{$k=1$} & \multicolumn{8}{c}{$k=2$} \\
    \cmidrule(lr){3-10} \cmidrule(lr){11-18}
    & & \multicolumn{4}{c|}{CFL = 0.5} & \multicolumn{4}{c||}{CFL = 2.0} & \multicolumn{4}{c|}{CFL = 0.5} & \multicolumn{4}{c}{CFL = 2.0} \\
    \cmidrule(lr){3-6} \cmidrule(lr){7-10} \cmidrule(lr){11-14} \cmidrule(lr){15-18}
    $\varepsilon$ & $N_x$ & $e^1$ & Ord & $e^2$ & Ord & $e^1$ & Ord & $e^2$ & Ord & $e^1$ & Ord & $e^2$ & Ord & $e^1$ & Ord & $e^2$ & Ord \\
    \midrule
    \multirow{3}{*}{$1$} 
    & 8  & 2.266e-02 & --- & 2.580e-02 & --- & 1.756e-02 & --- & 1.892e-02 & --- & 3.195e-03 & --- & 3.351e-03 & --- & 1.876e-03 & --- & 2.113e-03 & --- \\
    & 16 & 6.807e-03 & 1.73 & 7.108e-03 & 1.86 & 5.052e-03 & 1.80 & 5.311e-03 & 1.83 & 5.173e-04 & 2.63 & 5.701e-04 & 2.56 & 3.193e-04 & 2.55 & 3.601e-04 & 2.55 \\
    & 32 & 1.781e-03 & 1.93 & 1.880e-03 & 1.92 & 1.382e-03 & 1.87 & 1.462e-03 & 1.86 & 5.841e-05 & 3.15 & 6.694e-05 & 3.09 & 3.305e-05 & 3.27 & 3.824e-05 & 3.24 \\
    \midrule
    \multirow{3}{*}{$10^{-2}$} 
    & 8  & 2.504e-02 & --- & 2.818e-02 & --- & 2.467e-02 & --- & 2.759e-02 & --- & 3.535e-03 & --- & 3.699e-03 & --- & 3.553e-03 & --- & 3.709e-03 & --- \\
    & 16 & 7.467e-03 & 1.75 & 7.759e-03 & 1.86 & 6.859e-03 & 1.85 & 7.150e-03 & 1.95 & 5.666e-04 & 2.64 & 6.406e-04 & 2.53 & 5.584e-04 & 2.67 & 6.382e-04 & 2.54 \\
    & 32 & 1.907e-03 & 1.97 & 2.016e-03 & 1.94 & 1.690e-03 & 2.02 & 1.789e-03 & 2.00 & 6.205e-05 & 3.19 & 7.327e-05 & 3.13 & 5.334e-05 & 3.39 & 6.329e-05 & 3.33 \\
    \midrule
    \multirow{3}{*}{$10^{-6}$} 
    & 8  & 2.730e-02 & --- & 3.046e-02 & --- & 2.857e-02 & --- & 3.203e-02 & --- & 3.966e-03 & --- & 3.997e-03 & --- & 4.224e-03 & --- & 4.316e-03 & --- \\
    & 16 & 7.218e-03 & 1.92 & 7.613e-03 & 2.00 & 8.207e-03 & 1.80 & 8.697e-03 & 1.88 & 7.361e-04 & 2.43 & 8.452e-04 & 2.24 & 8.254e-04 & 2.36 & 9.440e-04 & 2.19 \\
    & 32 & 1.736e-03 & 2.06 & 1.817e-03 & 2.07 & 2.102e-03 & 1.97 & 2.218e-03 & 1.97 & 1.347e-04 & 2.45 & 2.109e-04 & 2.00 & 1.064e-04 & 2.96 & 1.426e-04 & 2.73 \\
    \bottomrule
    \end{tabular}}
    \caption{\small Accuracy results for \texttt{ARS443} in \textbf{Test I}.}
    \label{table:ars443}
\end{table}
In Table~\ref{table:ars443}, we present the relative errors and convergence orders of the scheme defined by 
$\text{Order}_{N_x} = (\log (e_{N_x/2}) - \log (e_{N_x}))/\log 2\,.$

From the table, we observe that the scheme exhibits consistent second order convergence for the case $k=1$. The convergence order is preserved across different $\varepsilon$ regimes, and is not affected by the use of a big CFL number. When $k=2$, the scheme achieves its expected third-order convergence when $\varepsilon=1$ and $\varepsilon=10^{-2}$, as shown in the right big column of Table~\ref{table:ars443}. However, we detect a reduction of order to $2$ when $\varepsilon=10^{-6}$. This loss of accuracy, however, is consistent with our observations in section 5.2.

\subsection*{Test II: AP property}
In this test, we analyze the AP property of our scheme by investigating the ability of the scheme to relax toward the equilibrium state. The AP error that measures the distance between the distribution function $f$ and the associated Maxwellian $M_f$ is defined as $\|f - M_f \|_{\ell^1}$. 
In order to observe a relaxation process toward the equilibrium state, we consider a non-equilibrium initial data
{\small
\begin{equation*}
    f (t=0, x, v) = \frac{\rho(x)}{4\pi T(x)} \exp \left(- \frac{(v - U(x))^2}{2T(x)}\right) + \frac{\rho(x)}{4\pi T(x)} \exp \left(- \frac{(v - W(x)))^2}{2T(x)}\right), 
\end{equation*}
}{\emergencystretch=2em Here, $\rho$ and $T$ are the same as in~\eqref{eq:macro_in_smooth_test}, and $U(x) = (5/4, -3/4)$, $W(x) = (-9/20, -3/4)$. Other settings are the same as in \textbf{Test I}, and we let $\varepsilon=10^{-2}, 10^{-4}, 10^{-6}$.\par}
 
We consider the methods \texttt{ARS443} and \texttt{DP2A242}, which are of type CK and type A, respectively. Theorem~\ref{AP property for schemes} tells us that the type-CK method \texttt{ARS443} leads to a weak AP property while the type-A method \texttt{DP2A242} is strongly AP. 
\begin{figure}
\centering
\includegraphics[width=0.8\textwidth]{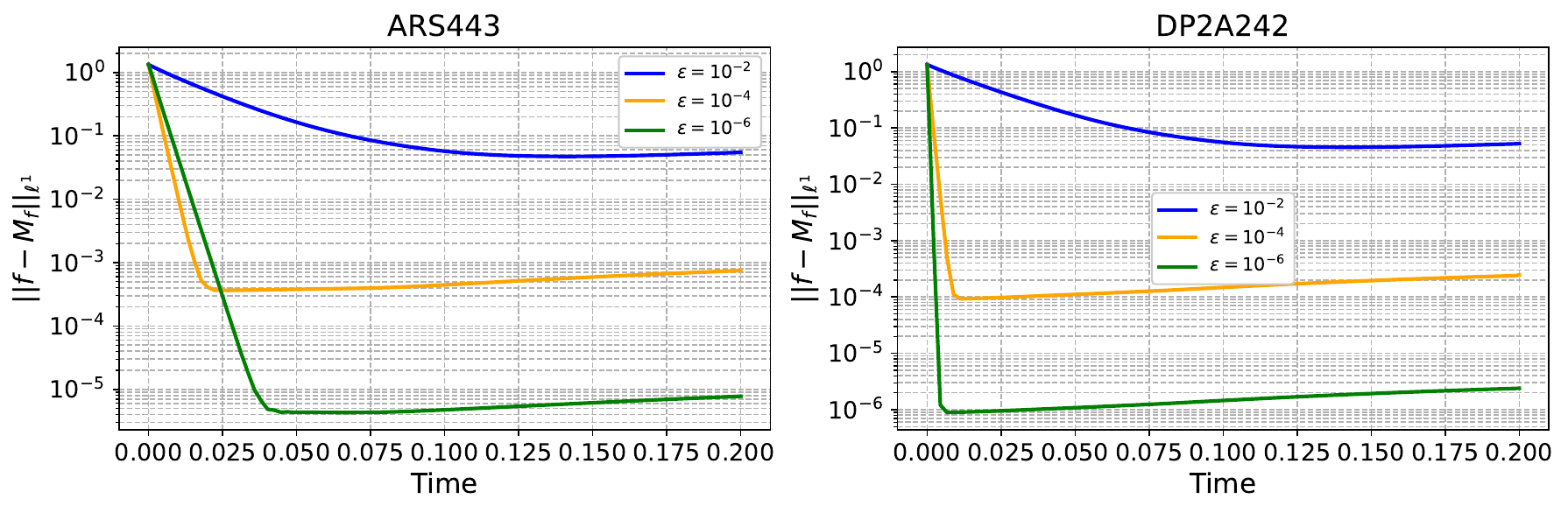}
    \caption{\small The time evolution of $||f - M_f||_{\ell^1}$ for \texttt{ARS443} and \texttt{DP2A242} in \textbf{Test II}.}
    \label{fig:aperror}
\end{figure}
In Figure~\ref{fig:aperror}, we find that for all $\varepsilon$, the distances between $f$ and its equilibrium $M_f$ diminish as time propagates and they reach the saturated values roughly at the same level as the corresponding $\varepsilon$. Secondly, the smaller $\varepsilon$ is, generally the faster $f$ reaches its equilibrium state. On the other hand, a stronger AP property is demonstrated in the type-A method \texttt{DP2A242}, as shown on the right side of Figure~\ref{fig:aperror}. In the cases where $\varepsilon=10^{-4}$ and $10^{-6}$, we observe that the errors $\|f - M_f \|_{\ell^1}$ rapidly drop and converge to their corresponding $\varepsilon$ values only after {\it one time step}. These observations are consistent with the conclusion in Theorem~\ref{AP property for schemes}. 

\subsection*{Test III: Sod shock problem}

In this test, we study the benchmark Riemann shock problem to show that our proposed schemes~\eqref{eq:IMEX-RK} can accurately capture the sharp discontinuities. The initial condition is given by
\[
\begin{aligned}
f(t=0, x, v) &= \frac{\rho(x)}{2\pi T(x)} \exp \left(- \frac{(v - U(x))^2}{2T(x)}\right),\\
(\rho, U, T) &= 
\begin{cases} 
(1, (0, 0), 1), & \text{if } 0 \leq x \leq 0.5, \\ 
\left( \frac{1}{8}, (0, 0), \frac{1}{4} \right), & \text{if } 0.5 < x \leq 1, 
\end{cases}
\end{aligned}
\]
Neumann boundary condition is assumed in the spatial direction. To mitigate spurious oscillations near discontinuities, we utilize a linear rescaling with LMPP limiter~\eqref{eq:LMPP} on the reconstructed polynomials. The numerical parameters are set as $N_x = 80$, $N_v = 32$ and $k = 2$. We use CFL values $0.5$ or $2$ and compare the results in different regimes at the time $t=0.2$, characterized by $\varepsilon = 10^{-2}$ (moderately stiff) and $\varepsilon = 10^{-8}$ (highly stiff). Reference solutions are obtained from the \texttt{AP-FD} method ~\cite{Filbet-Jin}, using $N_x = 200$ and $\Delta t = 3 \times 10^{-4}$.

\begin{figure}[!h]
    \centering
    \begin{subfigure}[b]{0.8\textwidth}
        \includegraphics[width=\textwidth]{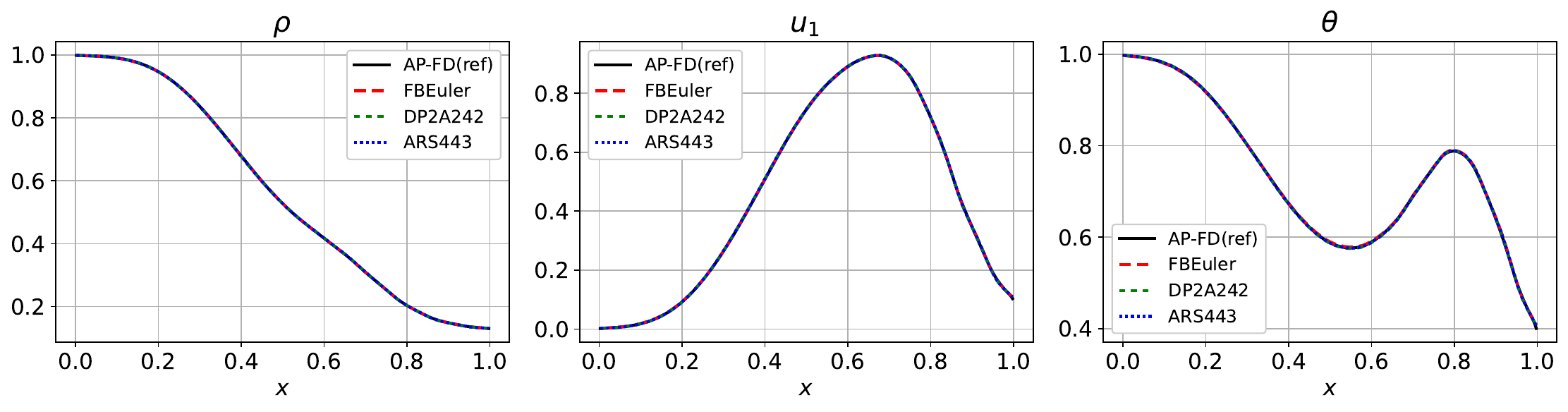}
        \caption{$\varepsilon = 10^{-2}$}
    \end{subfigure}
    \vfill
    \begin{subfigure}[b]{0.8\textwidth}
        \includegraphics[width=\textwidth]{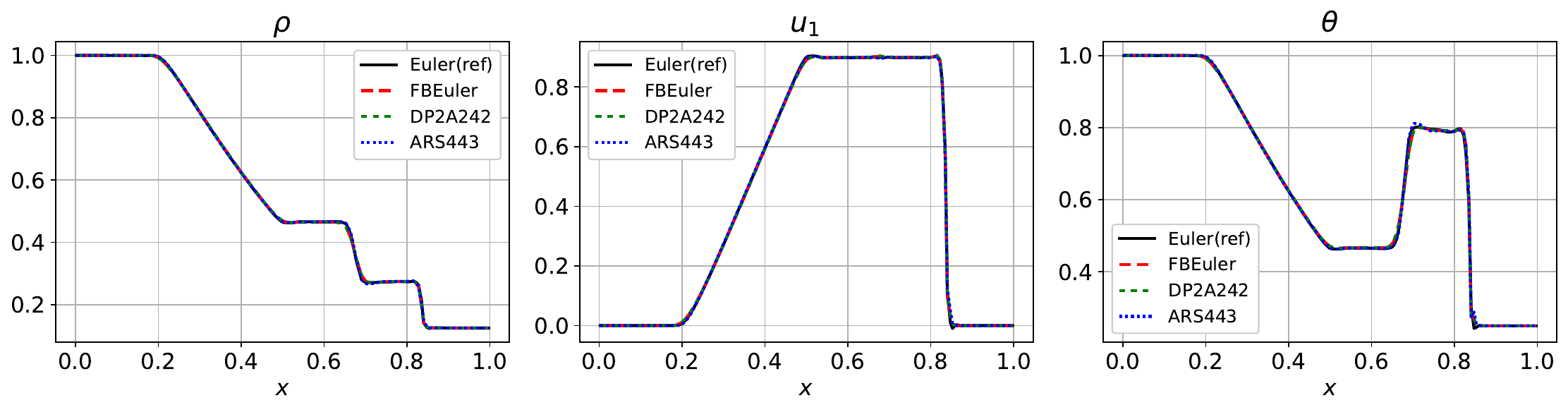}
        \caption{$\varepsilon = 10^{-8}$}
    \end{subfigure}
    \caption{\small Density, the first component of
mean velocity and temperature of \textbf{Test III}, with CFL$=0.5$.}
    \label{fig:shock0.5}
\end{figure}

\begin{figure}
    \centering
    \begin{subfigure}[b]{0.8\textwidth}
        \includegraphics[width=\textwidth]{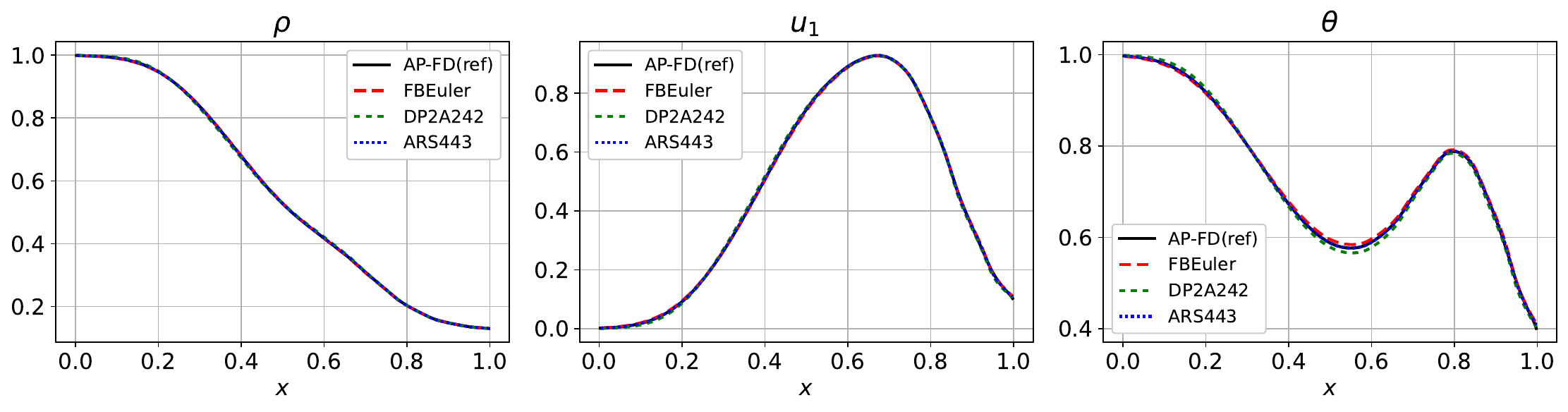}
        \caption{$\varepsilon = 10^{-2}$}
    \end{subfigure}
    \vfill
   \begin{subfigure}[b]{0.8\textwidth}
        \includegraphics[width=\textwidth]{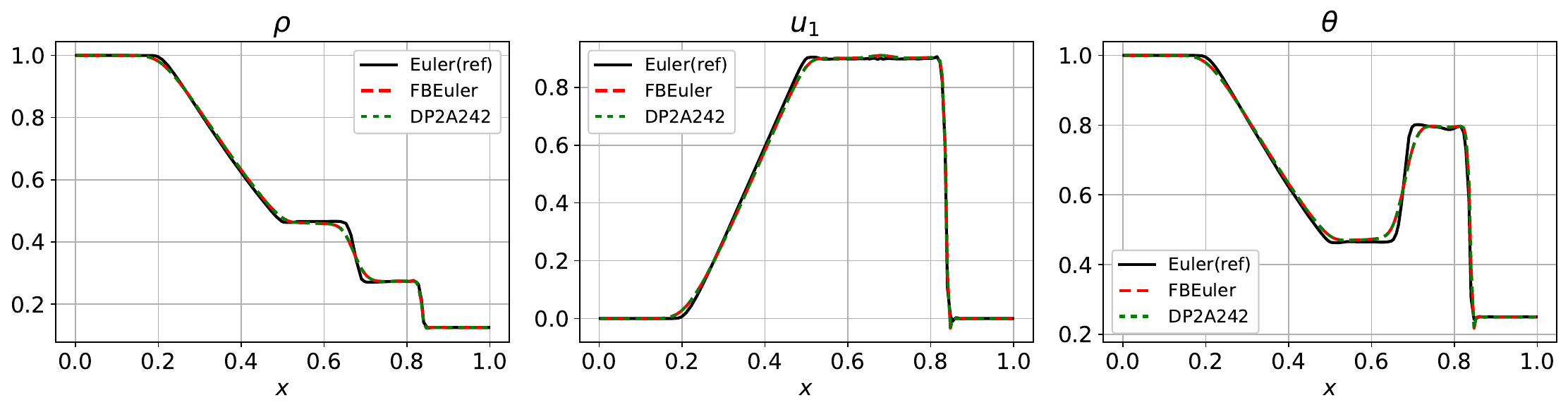}
        \caption{$\varepsilon = 10^{-8}$}
    \end{subfigure}
    \caption{\small Density, the first component of
mean velocity and temperature of \textbf{Test III}, with CFL$=2$. } 
    \label{fig:shock2}
\end{figure}

In Figure~\ref{fig:shock0.5}, under $\text{CFL} = 0.5$ the numerical solutions computed by using all three Butcher tableaus show good agreements with the reference solutions. 
Moreover, in Figure~\ref{fig:shock2}, under $\text{CFL}=2$ and the challenging case of $\varepsilon=10^{-8}$, the schemes \texttt{FBEuler} and \texttt{DP2A242}
remain stable and relatively accurate--though a slight mismatch is observed--under larger time steps. This attributes to the CFL-independent property of the scheme, which is considered as a major advantage of SL methods.

Applying higher-order schemes such as \texttt{ARS443} to the Riemann shock problem presents additional challenges. Beyond the AP property, the positivity-preserving property is crucial for ensuring stability.
While the LMPP limiter is applied in the spatial direction to prevent overshoots, the temporal reconstruction introduces further requirements. As established in Theorem~\ref{positivity constraints for RK}, \texttt{ARS443} does not strictly satisfy the sufficient constraints on the coefficients to ensure positivity. Our numerical experiments confirm that at $\text{CFL}=2$, the \texttt{ARS443} scheme works for relatively large $\varepsilon$ but becomes unstable for small $\varepsilon$ due to positivity loss, whereas it remains stable at $\text{CFL}=0.5$ for all $\varepsilon$ values.

\subsection*{Test IV: A mixing regime test}

In this section, we test the performance of the proposed method in a case with a spatially varying $\varepsilon(x)$ given by
{\small
\begin{equation}
    \label{eq:SpatiallyVaryingEpsilon}
    \varepsilon(x) = \varepsilon_0 + \frac{1}{2}\left( \tanh(1-10(x-0.5)) + \tanh(1+10(x-0.5)) \right), 
\end{equation}
}where \(\varepsilon_0 = 1 \times 10^{-6}\). This problem varies from the fluid to the kinetic regime. 
Assume that the initial data is far away from the local equilibrium and defined by
{\small
\begin{equation*}
    f(t=0, x, v) = \frac{\rho(x)}{4\pi T(x)} \exp \left(- \frac{(v - U(x))^2}{2T(x)}\right) + \frac{\rho(x)}{4\pi T(x)} \exp \left(- \frac{(v + U(x))^2}{2T(x)}\right), 
\end{equation*}
}where 
$\rho(x) = \frac{2 + \sin (2\pi x)}{3}, \; U(x) = (\cos (2\pi x), 0), \; T(x) = \frac{3 + \cos (2\pi x)}{4}. $
Periodic boundary condition is assumed in the spatial dimension. 
We study and compare the performance of \texttt{FBEuler}, \texttt{DP2A242}, and \texttt{ARS443} methods. 
Similar to \textbf{Test III}, reference solutions are obtained from the \texttt{AP-FD} method ~\cite{Filbet-Jin}, using $N_x = 1000$ and $\Delta t = 6 \times 10^{-5}$. 
\begin{figure}
    \centering
    
    \begin{subfigure}[b]{0.8\textwidth}
        \includegraphics[width=\textwidth]{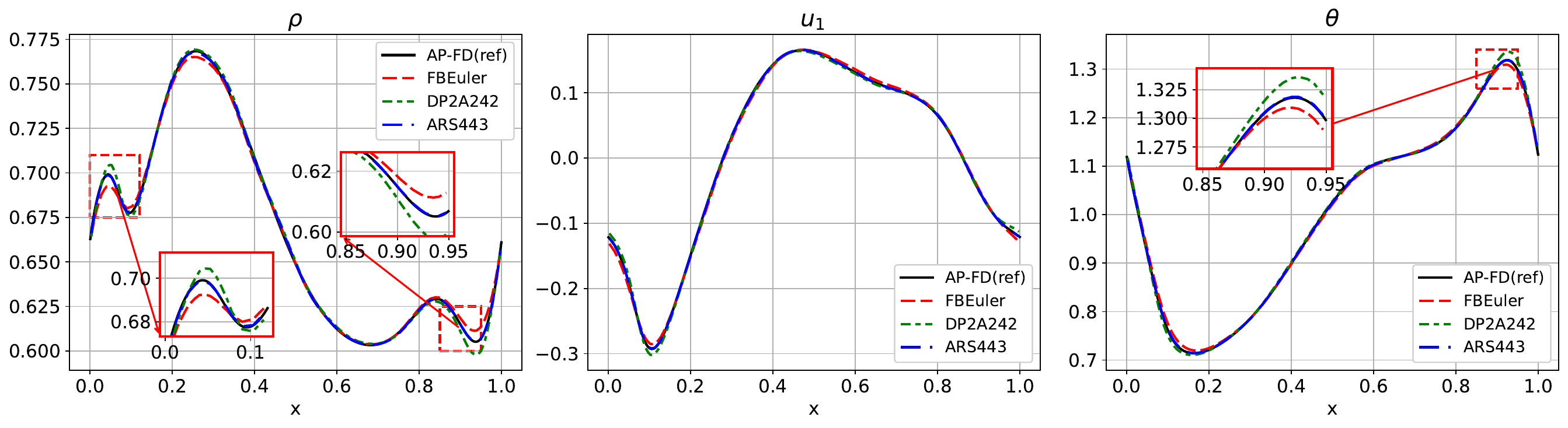}
    \end{subfigure}
    
    \caption{\small Density, the first component of
mean velocity and temperature of \textbf{Test IV} at $t=0.3$, with CFL = 0.5. 
     }
    \label{fig:mixedRegime0.5}
\end{figure}
\begin{figure}
    \centering  
    
    \begin{subfigure}[b]{0.8\textwidth}
        \includegraphics[width=\textwidth]{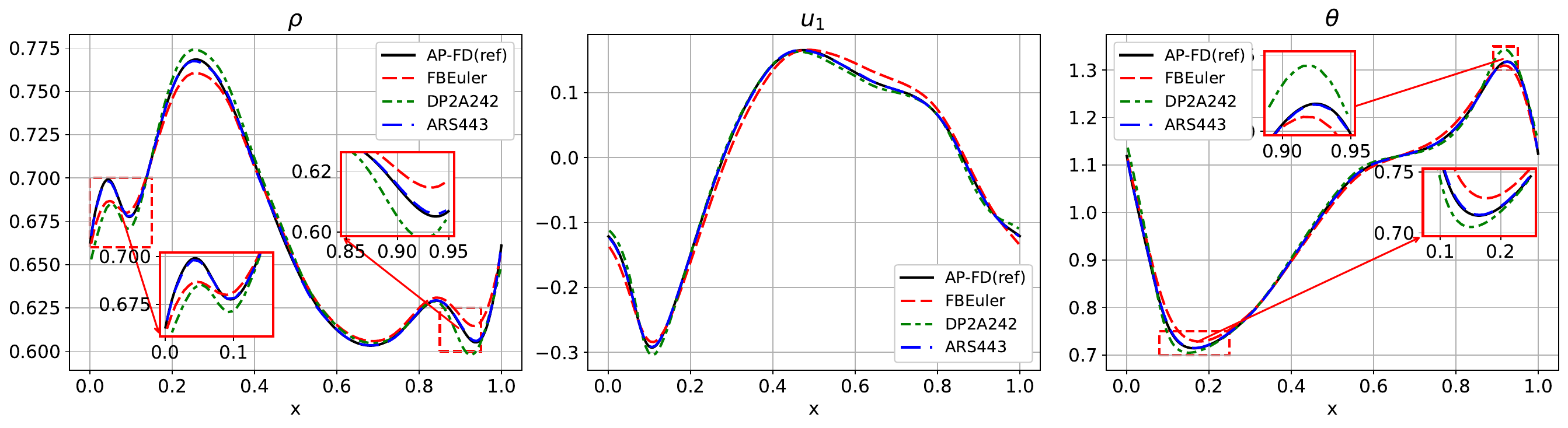}
    \end{subfigure}
    
    \caption{\small Density, the first component of
mean velocity and temperature of \textbf{Test IV} at $t=0.3$, with CFL = 2. }
    \label{fig:mixedRegime2}
\end{figure}
In Figure~\ref{fig:mixedRegime0.5} and Figure~\ref{fig:mixedRegime2}, we compare density, the first component of mean velocity and temperature with the reference solutions at time $t=0.3$. When $\text{CFL}=0.5$, all three numerical solvers demonstrate good agreements with their reference solutions. 
The first-order \texttt{FBEuler} shows reasonable accuracy, while the second-order \texttt{DP2A242} and the third-order \texttt{ARS443} schemes produce increasingly accurate results. 
When the CFL number increases to 2, accuracy of the solution continues to enhance as the scheme order increases. However, some discrepancies of the numerical solution from the reference solution start to emerge, particularly among lower-order methods. Both \texttt{FBEuler} and \texttt{DP2A242} methods show substantial accuracy degradation under the more demanding time step, while macroscopic quantities solved by the \texttt{ARS443} method maintain close agreements with the reference solution, demonstrating its robustness even when CFL is relatively large.

\section{Conclusion}
\label{sec:conclusion}
We develop and analyze asymptotic-preserving semi-Lagrangian discontinuous Galerkin schemes for the Boltzmann equation, together with a derivation of the limiting system as $\varepsilon \to 0$. Our analysis reveals that existing IMEX-RK schemes retain their temporal accuracy for fixed $\varepsilon$ in the SL setting, while order reduction is still present as $\varepsilon \to 0$ and additional conditions on the RK coefficients need to be imposed to ensure higher order accuracy in the limiting scheme. While semi-Lagrangian methods eliminate CFL restrictions, positivity preservation may impose alternative constraints. The hypocoercivity-based stability analysis provides foundations for long-time behavior, and numerical experiments validate the scheme's advantages and limitations across diverse regimes. Future work includes extending stability analysis to the nonlinear equation, designing IMEX-RK tableaux satisfying high-order, asymptotic accuracy, and positivity constraints simultaneously, and developing adaptive strategies for multi-scale features.

\section*{Acknowledgments}
The authors thank Jiaojiao Guo, a PhD student in the Department of Mathematics, Faculty of Arts and Sciences, Beijing Normal University, for her helpful suggestions on revising the manuscript.

\appendix

\section{DG Spatial Discretization for~\eqref{eq:one_step_discretization}}
\label{sec:DG}

\begin{figure}[htbp]
    \centering
    \resizebox{0.95\textwidth}{!}{%
    \begin{tikzpicture}[
        x=3.2cm,
        y=2.6cm,
        every node/.style={font=\small},
        axis/.style={draw=axiscolor, line width=0.8pt},
        characteristic/.style={
            draw=footblue,
            line width=1pt,
            dashed,
            -{Stealth[length=2.2mm,width=1.6mm]}
        },
        dimension/.style={
            line width=0.7pt,
            {Stealth[length=1.5mm]}-{Stealth[length=1.5mm]}
        }
    ]

        \definecolor{axiscolor}{HTML}{1D2A3A}
        \definecolor{footblue}{HTML}{2463C4}
        \definecolor{cellorange}{HTML}{EF9325}
        \definecolor{cellteal}{HTML}{11968F}
        \definecolor{labelred}{HTML}{C92B3F}
        \definecolor{dimgray}{HTML}{8590A0}

        \pgfmathsetmacro{\alphaVal}{0.3}
        \pgfmathsetmacro{\footL}{-0.5+\alphaVal}
        \pgfmathsetmacro{\footR}{ 0.5+\alphaVal}
        \pgfmathsetmacro{\footMid}{(\footL+\footR)/2}

        \foreach \x in {-1.5,-0.5,0.5,1.5} {
            \draw[
                densely dotted,
                draw=dimgray!35,
                line width=0.5pt
            ] (\x,0) -- (\x,0.7);
        }

        \path[
            fill=footblue!8,
            draw=footblue,
            line width=0.7pt
        ] (0.5,0.67) rectangle (1.5,0.73);

        \draw[axis] (-1.7,0.7) -- (1.65,0.7)
            node[right=7pt, text=black] {$t^{n+1}$};

        \draw[axis] (-1.7,0) -- (1.65,0)
            node[right=7pt, text=black] {$t^n$};

        \draw[
            draw=cellorange,
            line width=3pt,
            line cap=butt
        ] (\footL,0) -- (0.5,0);

        \draw[
            draw=cellteal,
            line width=3pt,
            line cap=butt
        ] (0.5,0) -- (\footR,0);

        \foreach \x in {-1.5,-0.5,0.5,1.5} {
            \draw[axis] (\x,0.735) -- (\x,0.665);
            \draw[axis] (\x,0.035) -- (\x,-0.035);
        }

        \node[above=6pt] at (-1.5,0.735) {$x_{j-5/2}$};
        \node[above=6pt] at (-0.5,0.735) {$x_{j-3/2}$};
        \node[above=6pt] at ( 0.5,0.735) {$x_{j-1/2}$};
        \node[above=6pt] at ( 1.5,0.735) {$x_{j+1/2}$};

        \node[above=5pt, text=footblue]
            at (1,0.735) {$I_j$};

        \draw[characteristic] (0.5,0.67) -- (\footL,0);
        \draw[characteristic] (1.5,0.67) -- (\footR,0);

        \draw[dimension, draw=dimgray]
            (-0.5,0.115) -- (\footL,0.115)
            node[midway, above=3pt, text=dimgray]
            {$\alpha\Delta x$};

        \node[below=3pt, text=labelred]
            at (-0.5,-0.075) {$x_{j^\ast-1/2}$};

        \node[below=3pt, text=labelred]
            at (0.5,-0.075) {$x_{j^\ast+1/2}$};

        \draw[
            densely dotted,
            draw=footblue,
            line width=0.6pt
        ] (\footL,-0.045) -- (\footL,-0.30);

        \draw[
            densely dotted,
            draw=footblue,
            line width=0.6pt
        ] (\footR,-0.045) -- (\footR,-0.30);

        \fill[footblue] (\footL,0) circle[radius=1.7pt];
        \fill[footblue] (\footR,0) circle[radius=1.7pt];

        \node[below=2pt, text=footblue]
            at (\footL,-0.33) {$x^*_{j-1/2}$};

        \node[below=2pt, text=footblue]
            at (\footR,-0.33) {$x^*_{j+1/2}$};

        \draw[dimension, draw=footblue]
            (\footL,-0.56) -- (\footR,-0.56);

        \node[below=3pt, text=footblue]
            at (\footMid,-0.63)
            {$I_j^\ast=
              [x^*_{j-1/2},x^*_{j+1/2}]$};

    \end{tikzpicture}%
    }
    \caption{Illustration of the upstream characteristic feet.}
    \label{fig:upstream_feet}
\end{figure}
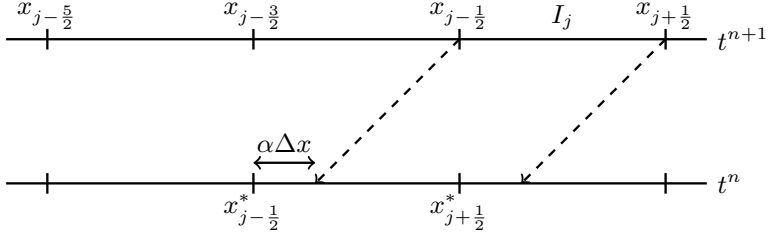

Let the test function be the Lagrange polynomials introduced in section~\ref{sec:nodalDG}, namely $\Psi(x) = \ell_{j, p_j}(x)$. We insert the representation~\eqref{eq:solution_representation} into~\eqref{eq:one_step_discretization}, and perform Gauss quadrature to evaluate the integral for terms at $t^{n+1}$. This leads to
$$
\begin{aligned}
    \Delta x\, \Theta^{n+1}_{j, p_j} \omega_{p_j} 
    &= \sum_{p} \Xi^n_{j^\ast\!,\, p} 
    \,\int_{x_{j^\ast - 1/2} + \alpha \Delta x}^{x_{j^\ast + 1/2}} 
    \ell_{j^\ast\!,\, p}(x) \ell_{j, p_j}(x + v \Delta t) \,\mathrm{d}x \\
    &\quad + \sum_{p} \Xi^n_{j^\ast + 1\!,\, p} 
    \,\int_{x_{j^\ast + 1/2}}^{x_{j^\ast + 1/2} + \alpha \Delta x} 
    \ell_{j^\ast + 1\!,\, p}(x) \ell_{j, p_j}(x + v \Delta t) \,\mathrm{d}x,
\end{aligned}
$$
where $\omega_{p_j}$ are the Gauss-Legendre weights on $[0, 1]$. The index $j^\ast$ and $\alpha$ satisfy $\displaystyle x_{j - 1/2} - v \Delta t = x_{j^{\ast} - 1/2} + \alpha \Delta x$, $\alpha\in [0,1)$, see Figure~\ref{fig:upstream_feet}. 
By a change of variables $\displaystyle x = x_{j^\ast - 1/2} + s \Delta x $, the equation above is transformed into
$$
\begin{aligned}
    \Delta x\, \Theta^{n+1}_{j, p_j} \omega_{p_j} 
    &= \Delta x \sum_{p} \Xi^n_{j^\ast\!,\, p} 
    \int_{s = \alpha}^1
    \tilde{\ell}_{p}(s) \tilde{\ell}_{p_j}(s - \alpha) \,\mathrm{d}s \\ & + \Delta x\sum_{p} \Xi^n_{j^\ast + 1\!,\, p} 
    \int_{s = 0}^\alpha
    \tilde{\ell}_{p}(s) \tilde{\ell}_{p_j}(s + 1 - \alpha) \,\mathrm{d}s,
\end{aligned}
$$
where $\tilde{\ell}_{p}(s)$ represents the Lagrange polynomials interpolating $k+1$ Gauss-Legendre points on $[0, 1]$. 
To evaluate the integrals, we use further changes of variables $s = \alpha + u(1 - \alpha)$ in the first integral, and $s = \alpha u$ in the second integral. Then
$$
\begin{aligned}
    \Delta x\, \Theta^{n+1}_{j, p_j} \omega_{p_j} 
    &= (1 - \alpha)\Delta x\sum_{p} \Xi^n_{j^\ast\!,\, p} 
    \int_{u = 0}^1
    \tilde{\ell}_{p}(\alpha + u(1 - \alpha)) \tilde{\ell}_{p_j}(u(1 - \alpha)) \,\mathrm{d}u \\
    &\quad + \alpha\Delta x \sum_{p} \Xi^n_{j^\ast + 1\!,\, p} 
    \int_{u = 0}^1
    \tilde{\ell}_{p}(\alpha u) \tilde{\ell}_{p_j}(\alpha(u - 1) + 1) \,\mathrm{d}u.
\end{aligned}
$$
These integrals can be computed by the Gauss-Legendre quadrature on $[0, 1]$. Applying the quadrature weights $\omega_q$ and nodes $u_q$, we obtain
$$
\begin{aligned}
    \Delta x\, \Theta^{n+1}_{j, p_j} \omega_{p_j} 
    &= (1 - \alpha)\Delta x \sum_{p}\sum_q \Xi^n_{j^\ast\!,\, p} 
    \, \omega_q \tilde{\ell}_{p}(\alpha + u_q (1 - \alpha)) \tilde{\ell}_{p_j}(u_q (1 - \alpha)) \\
    &\quad + \alpha \Delta x \sum_{p}\sum_q \Xi^n_{j^\ast + 1\!,\, p} 
    \, \omega_q \tilde{\ell}_{p}(\alpha u_q) \tilde{\ell}_{p_j}(\alpha(u_q - 1) + 1).
\end{aligned}
$$
Finally, writing in terms of the original variables, one obtains
$$
\begin{aligned}
    f^{n+1}_{j, \cdot} &=  \Delta t \frac{\beta^{n+1}_{j, \cdot}}{\varepsilon} (\mathcal{M}^{n+1}_{j, \cdot} - f^{n+1}_{j, \cdot}) + A(\alpha) \left[ f^n_{j^\ast, \cdot} + \frac{\Delta t}{\varepsilon} (Q^n_{j^\ast, \cdot} - \beta^n_{j^\ast, \cdot} (\mathcal{M}^n_{j^\ast, \cdot} - f^n_{j^\ast, \cdot})) \right] \\
    &\quad + B(\alpha) \left[ f^n_{j^\ast + 1, \cdot} + \frac{\Delta t}{\varepsilon} (Q^n_{j^\ast + 1, \cdot} - \beta^n_{j^\ast + 1, \cdot} (\mathcal{M}^n_{j^\ast + 1, \cdot} - f^n_{j^\ast + 1, \cdot})) \right],
\end{aligned}
$$
where $Q^n = Q(f^n)$ and $\mathcal{M}^n = \mathcal{M}(f^n)$.

\section{Definitions} 
\label{sec:definitions}
In this section, we give the definitions of various IMEX-RK methods as in~\cite{Dimarco-Pareschi}.
\begin{definition}\label{def A-CK schemes}
    We call an IMEX-RK method of type A if the matrix $\tilde{A} \in \mathbb{R}^{s \times s}$ is invertible. We call an IMEX-RK method of type CK if the matrix $\tilde{A}$ can be written as 
    \begin{equation*}
        \tilde{A} = \begin{pmatrix}
            0 & 0 \\
            \tilde{a} & \tilde{\hat{A}}
        \end{pmatrix}
    \end{equation*}
    where $\tilde{a} \in \mathbb{R}^{(s-1)}$ and the submatrix $\tilde{\hat{A}} \in \mathbb{R}^{(s-1) \times (s-1)}$ is invertible. In the special case $\tilde{a}=0, \tilde{b}_1 =0$ the scheme is said to be of type ARS. When the scheme is of type CK, we write 
    \begin{equation*}
        A = \begin{pmatrix}
            0 & 0 \\
            a & \hat{A}
        \end{pmatrix}
\end{equation*}
to make the notation of $A$ consistent with that of $\tilde{A}$. 
\end{definition}
\begin{definition}\label{def GSA}
    We call an IMEX-RK method globally stiffly accurate (GSA) if $a_{si}=b_i$ and $\tilde{a_{si}} = \tilde{b}_{i}$ $\forall 1 \leq i \leq s$. That is, the numerical solution is equal to the last stage value. 
\end{definition}

\begin{definition}
\label{def:matrices}

{\small
\begin{itemize}
    \item For any vector $w=(w_1,\ldots,w_n)^{T} \in \mathbb{R}^{n}$:
    \begin{itemize}
        \item $w_{(i)} \in \mathbb{R}^{i}$: vector consisting of the first $i$ entries of $w$
        \item $\hat{w} \in \mathbb{R}^{n-1}$: vector with the first entry of $w$ deleted
    \end{itemize}
    \item For any matrix $C \in \mathbb{R}^{n \times n}$:
    \begin{itemize}
        \item $C_{(i)}$: submatrix of $C$ consisting of the first $i$ rows and $i$ columns
        \item $C_{i}$: the $i$-th row of $C_{(i)}$
        \item $\bar{C}_{i} = (c_{i1},c_{i2},\ldots,c_{i(i-1)}) \in \mathbb{R}^{(i-1)}$: the $i$-th row ($i\geq2$) of $C$, with only the first $(i-1)$ entries
    \end{itemize}
    \item For invertible lower-triangular matrix $C$:
    \begin{itemize}
        \item $(C^{-1})_{(i)} = (C_{(i)})^{-1}$ for all $1 \leq i \leq n$
        \item $C^{-1}_{(i)}$: represents either of the two equivalent quantities above
    \end{itemize}
    \item $e=(1,1,\ldots,1)^{T} \in \mathbb{R}^{s}$, $F=(f^{(1)},f^{(2)},\ldots,f^{(s)})^{T}$: vector of functions
    \item For any $\mathcal{F}=(f^{j_1},f^{j_2},\ldots,f^{j_k})^{T}$ with $1 \leq j_1,j_2,\ldots,j_k \leq s$:
    \begin{align*}
        \mathcal{S}^{i}[\mathcal{F}] &:=(\mathcal{S}_{i,j_1}[f^{(j_1)}],\mathcal{S}_{i,j_2}[f^{(j_2)}],\ldots,\mathcal{S}_{i,j_k}[f^{(j_k)}])^{T}, \; \;
        \tilde{\mathcal{S}}^{i}[\mathcal{F}] \text{ defined similarly}
    \end{align*}
\end{itemize}
}    
\end{definition}

\begin{definition}
\label{def:coeff_order_condition_moments}
The coefficients for Theorem~\ref{thm:order_condition_moments} are given by
{\small
\begin{align*}
& \tilde{b}_{kj} = \frac{\tilde{a}_{kj}}{\tilde{a}_{jj}} - \sum_{l=j+1}^{k-1} \frac{\tilde{a}_{kl} \tilde{b}_{lj}}{\tilde{a}_{ll}}, \quad k > j \geq 1, \quad C_k = \tilde{c}_k, \quad D_k = \sum_{j=1}^{k-1} \tilde{b}_{kj} (D_j + (\tilde{c}_k - \tilde{c}_j)\tilde{c}_j), \\
& B_k = \left( 1 - \sum_{j=1}^{k-1} \tilde{b}_{kj} \right) \tilde{c}_k^2 + \sum_{j=1}^{k-1} \tilde{b}_{kj} (B_j + (\tilde{c}_k - \tilde{c}_j)^2), \quad G_k = \sum_{j=1}^{k-1} \tilde{b}_{kj} (G_j + \frac{1}{2} (\tilde{c}_k - \tilde{c}_j)^2 \tilde{c}_j), \\ & H_k = \sum_{j=1}^{k-1} \tilde{b}_{kj} (H_j + (\tilde{c}_k - \tilde{c}_j) D_j), \quad B_k^* = \sum_{j=1}^{k-1} \tilde{b}_{kj} (B_j^* + (\tilde{c}_k - \tilde{c}_j) B_j), \\ & 
B_k^{**} = \sum_{j=1}^{k-1} \tilde{b}_{kj} (B_j^{**} + (\tilde{c}_k - \tilde{c}_j)^2 \tilde{c}_j), \\ & B_k^{***} = \left( 1 - \sum_{j=1}^{k-1} \tilde{b}_{kj} \right) \tilde{c}_k^3 + \sum_{j=1}^{k-1} \tilde{b}_{kj} (B_j^{***} + (\tilde{c}_k - \tilde{c}_j)^3).
\end{align*}
}

\end{definition}

\section{Local-maximum-principle-preserving (LMPP) limiter} 
\label{sec:LMPP}

For problems with shocks, in order to control spurious oscillations near the spatial discontinuities, we apply the LMPP limiter used in~\cite{Ding-Qiu-Shu-2024, Zhang-Shu-2010}. 
Consider the cell $I_i$, since numerical solutions are represented as combinations of Lagrangian polynomials, we let $p(x)$ be such a polynomial in $I_i$, and denote $\bar{p}$ as its average over the cell. The modification consists in performing a linear rescaling to $p(x)$ as
\(
    \tilde{p}(x) := \theta (p(x) - \bar{p}) + \bar{p}, 
\)
with the limiter defined by
\begin{equation}
\label{eq:LMPP}
\theta = \min \left\{ \left|\frac{M^n - \bar{p}}{M^{n+1} - \bar{p}}\right|, \, \left|\frac{m^n - \bar{p}}{m^{n+1} - \bar{p}}\right|, 1 \right\},
\end{equation}
where $M^{n+1} = \max_{x\in I_i} p(x)$, $m^{n+1} = \min_{x\in I_i} p(x)$, $M^{n} = \max_{x\in I^\star_i} p(x)$, and $m^{n} = \min_{x\in I^\star_i} p(x)$.
$I^\star_i$ denotes the union of all background Eulerian cells that cover the cell obtained by back-tracing $I_i$. For example, in Figure~\ref{fig:upstream_feet}, $I^\star_i = [x_{j-\frac{1}{2}}^*, x_{j+\frac{1}{2}}^*]$ when $\alpha = 0$, $I^\star_i = [x_{j-\frac{1}{2}}^*, x_{j+\frac{3}{2}}^*]$ when $\alpha > 0$.

\section{IMEX-RK Butcher tableaux}
\label{sec:Butcher}

In this section, we present several useful Butcher tableaux for~\eqref{eq:IMEX-RK}. 
The first-order forward-backward Euler method corresponds to the tableaux
\begin{equation}
    \label{eq:IMEX-RK-1st}
    \begin{array}{cc}
        \begin{array}{c|cc}
            0 & 0 & 0 \\
            1 & 1 & 0 \\
            \hline
            & 1 & 0
        \end{array} & \hspace{5pt}
        \begin{array}{c|cc}
            0 & 0 & 0 \\
            1 & 0 & 1 \\
            \hline
            & 0 & 1
        \end{array}
    \end{array}
\end{equation}
The second-order $4$-stage DP2A242 method~\cite{Dimarco-Pareschi} with $\gamma = 2$ corresponds to the tableaux
\begin{equation}
    \label{eq:IMEX-RK-DP2A242}
\begin{array}{cc}
    \begin{array}{c|cccc}
        0 & 0 & 0 & 0 & 0 \\
        0 & 0 & 0 & 0 & 0 \\
        1 & 0 & 1 & 0 & 0 \\
        1 & 0 & 1/2 & 1/2 & 0 \\
        \hline
        & 0 & 1/2 & 1/2 & 0
    \end{array} & \hspace{5pt}
    \begin{array}{c|cccc}
        2 & 2 & 0 & 0 & 0 \\
        0 & -2 & 2 & 0 & 0 \\
        1 & 0 & -1 & 2 & 0 \\
        1 & 0 & 1/2 & 1/2-2 & 2 \\
        \hline
        & 0 & 1/2 & 1/2-2 & 2
    \end{array}
\end{array}
\end{equation}
The third order $5$-stage GSA ARS443 method~\cite{Ascher-Ruuth-Spiteri-1997, Xiong-Jang-Li-Qiu-2015} corresponds to the tableaux
\begin{equation}
    \label{eq:IMEX-RK-ARS443}
\begin{array}{cc}
    \begin{array}{c|ccccc}
        0 & 0 & & & & \\
        1/2 & 1/2 & 0 & & & \\
        2/3 & 11/18 & 1/18 & 0 & & \\
        1/2 & 5/6 & -5/6 & 1/2 & 0 & \\
        1 & 1/4 & 7/4 & 3/4 & -7/4 & 0 \\
        \hline
        & 1/4 & 7/4 & 3/4 & -7/4 & 0
    \end{array} & \hspace{5pt}
    \begin{array}{c|ccccc}
        0 & 0 & & & & \\
        1/2 & 0 & 1/2 & & & \\
        2/3 & 0 & 1/6 & 1/2 & & \\
        1/2 & 0 & -1/2 & 1/2 & 1/2 & \\
        1 & 0 & 3/2 & -3/2 & 1/2 & 1/2 \\
        \hline
        & 0 & 3/2 & -3/2 & 1/2 & 1/2
    \end{array}
\end{array}
\end{equation}

\section{Shu-Osher form of~\eqref{eq:IMEX-RK}}
\label{sec:ShuOsher}
We introduce the Shu-Osher form~\cite{Shu-Osher-1988} of the scheme~\eqref{eq:IMEX-RK}. This form rewrites the original scheme in an equivalent way, using coefficients computed from a recursive formula.

\begin{theorem}[Shu-Osher form of~\eqref{eq:IMEX-RK}]
\label{Shu-Osher form for IMEX-RK}
Adopting the convention that all quantities with non-positive subscripts are treated as zero, we have
\begin{itemize}
    \item If~\eqref{eq:IMEX-RK} is of \textbf{type A}, then it can be written as
  \small{ \begin{equation}\label{formula for Shu-Osher form for IMEX-RK type A}
        \begin{aligned}
            f^{(i)} =& \big(1-\bar{\tilde{A}}_{i}\tilde{A}^{-1}_{(i-1)}e_{(i)}\big)\tilde{\mathcal{S}}_{i,0}[f^{n}] + \big(\bar{\tilde{A}}_{i}\tilde{A}^{-1}_{(i-1)}\big)\tilde{\mathcal{S}}^{i}[F_{(i-1)}]+ \frac{\Delta t}{\varepsilon}\tilde{a}_{ii}Q_{P}(f^{(i)})\\
            &+ \Delta t E^{i}[\mathcal{S}^{i},\tilde{\mathcal{S}}^{i},\frac{1}{\varepsilon} G_{P}(F_{(i-1)})], \quad i=1, \dots, s,
        \end{aligned}
   \end{equation}}where $E^{i}[\mathcal{S}^{i},\tilde{\mathcal{S}}^{i}, G_{P}(F_{(i-1)})]$ is a linear combination of terms involving $\mathcal{S}_{i,j},\tilde{\mathcal{S}}_{i,j}$ and $G_{P}(f^{(j)})$ for $1 \leq j \leq i-1$. 
   
   \item If~\eqref{eq:IMEX-RK} is of \textbf{type CK} and $c_i=\tilde{c}_i ~\forall 1 \leq i \leq s$, then it can be written as
   {\small
\begin{equation}\label{formula for Shu-Osher form for IMEX-RK type CK}
    \begin{aligned}
        f^{(i)} =& \big(1-B_i\hat{e}_{(i-1)}\big)\mathcal{S}_{i,0}[f^{n}] + B_i\mathcal{S}^{i}[\hat{F}_{(i-1)}]+ \frac{\Delta t}{\varepsilon}\tilde{\hat{a}}_{(i-1)(i-1)}Q_{P}(f^{(i)})\\
        &+ \Delta t D_i\mathcal{S}^{i}[\frac{1}{\varepsilon} G_{P}(\hat{F}_{(i-1)})]+ \Delta t d_i\mathcal{S}_{i,0}[ \frac{1}{\varepsilon} G_{P}(f^{n})], \; i=1,\dots,s,
    \end{aligned}
\end{equation}
}where
$
    B_i = \bar{\tilde{\hat{A}}}_{i-1}\tilde{\hat{A}}^{-1}_{(i-2)},\; D_i = \bar{\hat{A}}_{i-1}-\bar{\tilde{\hat{A}}}_{i-1}\tilde{\hat{A}}^{-1}_{(i-2)}\hat{A}_{(i-2)},\; d_i = a_{i-1} - B_i a_{(i-2)}.
$
\end{itemize}

\begin{proof}
We present an induction argument to prove the Theorem. When $i=1$, $\tilde{A}_{(i-1)}=0$ and $f^{(1)} = \mathcal{S}_{1,0}[f^n] + \frac{\Delta t}{\varepsilon}\tilde{a}_{11}Q_{P}(f^{(1)})$ clearly holds by \eqref{eq:IMEX-RK}. For $i \geq 2$, we suppose \eqref{formula for Shu-Osher form for IMEX-RK type A} holds for all $j < i$, then we can derive 
{\small
\begin{equation}\label{Shu-Osher formula for Q_P(fi)}
    \begin{aligned}
        Q_{P}(f^{(j)}) =& \frac{\varepsilon}{\tilde{a}_{jj} \Delta t}\Big( f^{(j)} - \big(1-\bar{\tilde{A}}_{j}\tilde{A}^{-1}_{(j-1)}e_{(j)}\big)\tilde{\mathcal{S}}_{j,0}[f^{n}] - \big(\bar{\tilde{A}}_{j}\tilde{A}^{-1}_{(j-1)}\big)\tilde{\mathcal{S}}^{j}[F_{(j-1)}] \Big) \\ & - \frac{1}{\tilde{a}_{jj}}E^{j}[\mathcal{S}^{j},\tilde{\mathcal{S}}^{j}, G_{P}(F_{(j-1)})] 
    \end{aligned}
\end{equation}
}for all$j < i$. Substituting \eqref{Shu-Osher formula for Q_P(fi)} into the original scheme \eqref{eq:IMEX-RK}, we have
\begin{equation}\label{Shu-Osher intermediate 1}
    \begin{aligned}
        f^{(i)} =& \tilde{\mathcal{S}}_{i,0} [f^n] + \frac{\Delta t}{\varepsilon}\sum_{j=1}^{i-1} a_{ij} \mathcal{S}_{i,j} [G_{P}(f^{(j)})] + \frac{\Delta t}{\varepsilon}\sum_{j=1}^{i-1} \tilde{a}_{ij} \tilde{\mathcal{S}}_{i,j} [Q_{P}(f^{(j)})] \\ &\quad + \frac{\Delta t}{\varepsilon}\tilde{a}_{ii}Q_{P}(f^{(i)}) \\
        =& \tilde{\mathcal{S}}_{i,0} [f^n] + \frac{\Delta t}{\varepsilon}\sum_{j=1}^{i-1} a_{ij} \mathcal{S}_{i,j} [G_{P}(f^{(j)})] + \sum_{j=1}^{i-1}\frac{\tilde{a}_{ij}}{\tilde{a}_{jj}}\tilde{\mathcal{S}}_{i,j}[f^{(j)}] \\
        &- \sum_{j=1}^{i-1}\frac{\tilde{a}_{ij}}{\tilde{a}_{jj}}\big(1-\bar{\tilde{A}}_{j}\tilde{A}^{-1}_{(j-1)}e_{(j)}\big)\tilde{\mathcal{S}}_{i,0}[f^{n}] - \sum_{j=1}^{i-1}\frac{\tilde{a}_{ij}}{\tilde{a}_{jj}}\big( \bar{\tilde{A}}_{j}\tilde{A}^{-1}_{(j-1)}\big)\tilde{\mathcal{S}}^{i}[F_{(j-1)}]  \\
        & - \frac{\Delta t}{\varepsilon}\sum_{j=1}^{i-1}\frac{\tilde{a}_{ij}}{\tilde{a}_{jj}}\tilde{S}_{i,j}[E^{j}[\mathcal{S}^{j},\tilde{\mathcal{S}}^{j}, G_{P}(F_{(j-1)})]] + \frac{\Delta t}{\varepsilon}\tilde{a}_{ii}Q_{P}(f^{(i)})
    \end{aligned}
\end{equation}
where we have applied the fact that the shifting operators $\tilde{\mathcal{S}}_{i,j}$ are linear and their actions are transitive in the sense that $\tilde{\mathcal{S}}_{i,j} \circ \tilde{\mathcal{S}}_{j,k} = \tilde{\mathcal{S}}_{i,k}$. 
Note that for type A schemes, $c_i \neq \tilde{c}_i$ in general. Consequently, $\mathcal{S}_{i,j} \neq \tilde{\mathcal{S}}_{i,j}$, and the transitivity property fails for compositions of the form $\mathcal{S}_{i,j} \circ \tilde{\mathcal{S}}_{j,k}$. Deriving an explicit formula for $E^i$ is both challenging and unnecessary, as it requires resolving the intricate interactions between $\mathcal{S}_{i,j}$ and $\tilde{\mathcal{S}}_{i,j}$. Therefore, we retain $E^i$ as an unresolved term in the Shu--Osher form.

Since $\tilde{A}_{(j)}\tilde{A}^{-1}_{(j)}=I$, we have $\bar{\tilde{A}}_{j}\tilde{A}^{-1}_{(j-1)} + \tilde{a}_{jj}\bar{\tilde{A}}^{-1}_{j}=0$, then 
{\footnotesize
\begin{equation}\label{Shu-Osher intermediate 2}
    \begin{aligned}
        &\quad 
        \sum_{j=1}^{i-1}\frac{\tilde{a}_{ij}}{\tilde{a}_{jj}}\tilde{\mathcal{S}}_{i,j}[f^{(j)}] - \sum_{j=1}^{i-1}\frac{\tilde{a}_{ij}}{\tilde{a}_{jj}}\big( \bar{\tilde{A}}_{j}\tilde{A}^{-1}_{(j-1)}\big)\tilde{\mathcal{S}}^{i}[F_{(j-1)}] 
        \\ & = \sum_{j=1}^{i-1}\frac{\tilde{a}_{ij}}{\tilde{a}_{jj}}\Big( \tilde{\mathcal{S}}_{i,j}[f^{(j)}]- \big( \bar{\tilde{A}}_{j}\tilde{A}^{-1}_{(j-1)}\big)\tilde{\mathcal{S}}^{i}[F_{(j-1)}]\Big)  =\sum_{j=1}^{i-1}\tilde{a}_{ij}\Big( \frac{1}{\tilde{a}_{jj}}\tilde{\mathcal{S}}_{i,j}[f^{(j)}] + \bar{\tilde{A}}^{-1}_{j}\tilde{\mathcal{S}}^{i}[F_{(j-1)}]\Big) 
        \\ & = \sum_{j=1}^{i-1}\tilde{a}_{ij}\tilde{A}^{-1}_{j}\tilde{\mathcal{S}}^{i}[F_{(j)}]  = \big(\bar{\tilde{A}}_{i}\tilde{A}^{-1}_{(i-1)}\big)\tilde{\mathcal{S}}^{i}[F_{(i-1)}]. 
    \end{aligned}
\end{equation}
}Similarly, we have 
{\footnotesize
\begin{equation}\label{Shu-Osher intermediate 3}
    \begin{aligned}
        &\quad \tilde{\mathcal{S}}_{i,0} [f^n] - \sum_{j=1}^{i-1}\frac{\tilde{a}_{ij}}{\tilde{a}_{jj}}\big(1-\bar{\tilde{A}}_{j}\tilde{A}^{-1}_{(j-1)}e_{(j-1)}\big)\tilde{\mathcal{S}}_{i,0}[f^{n}] 
        \\ & = \Big( 1- \sum_{j=1}^{i-1} \tilde{a}_{ij}\big( \frac{1}{\tilde{a}_{jj}} +\bar{\tilde{A}}_{j}^{-1}e_{(j-1)}\big)\Big)\tilde{\mathcal{S}}_{i,0}[f^{n}] = \Big( 1- \sum_{j=1}^{i-1} \tilde{a}_{ij}\tilde{A}^{-1}_{j}e_{(j)}\Big)\tilde{\mathcal{S}}_{i,0}[f^{n}] \\ &
        = \big(1-\bar{\tilde{A}}_{i}\tilde{A}^{-1}_{(i-1)}e_{(i)}\big)\tilde{\mathcal{S}}_{i,0}[f^{n}]. 
    \end{aligned}
\end{equation}
}Finally, let 
{\footnotesize
\begin{equation}\label{Shu-Osher intermediate 4}
    E^{i}[\mathcal{S}^{i},\tilde{\mathcal{S}}^{i}, F_{(i-1)}] = \sum_{j=1}^{i-1} a_{ij} \mathcal{S}_{i,j} [G_{P}(f^{(j)})]- \sum_{j=1}^{i-1}\frac{\tilde{a}_{ij}}{\tilde{a}_{jj}}\tilde{S}_{i,j}[E^{j}[\mathcal{S}^{j},\tilde{\mathcal{S}}^{j}, G_{P}(F_{(j-1)})]]
\end{equation}
}Substituting \eqref{Shu-Osher intermediate 2}, \eqref{Shu-Osher intermediate 3} \eqref{Shu-Osher intermediate 4} into \eqref{Shu-Osher intermediate 1}, we get precisely \eqref{formula for Shu-Osher form for IMEX-RK type A} for $i$. 

If the scheme is of type CK, we can make a similar induction argument to derive \eqref{formula for Shu-Osher form for IMEX-RK type CK}. The base case for $i=1,2$ can be easily checked from \eqref{eq:IMEX-RK}. The inductive step is the same as the argument for type A scheme, so we skip the details here. 
\end{proof}
\end{theorem}

\bibliographystyle{amsplain}
\bibliography{ref}
\end{document}